\documentclass[reqno]{amsart}

\usepackage{amssymb}

\theoremstyle{definition}
\newtheorem{definition}{Definition}[section]
\theoremstyle{plain}
\newtheorem{theorem}[definition]{Theorem}
\newtheorem{lemma}[definition]{Lemma}
\newtheorem{corollary}[definition]{Corollary}
\newtheorem{proposition}[definition]{Proposition}

\theoremstyle{remark}

\numberwithin{equation}{section}

\newcommand{\bsn}{\bigskip}

\newcommand{\F}{\ensuremath{{\mathcal F}}}
\newcommand{\Fl}{\F_\lambda}
\newcommand{\M}{\ensuremath{{\mathcal M}}}
\newcommand{\HM}{{}_H\M}
\newcommand{\HMH}{{}_H\M^H_H}
\renewcommand{\L}{\ensuremath{{\mathcal L}}}
\newcommand{\McoH}{M^{coH}}
\newcommand{\dH}{\hat{H}}

\newcommand{\RR}{{\rm \mathbf R}}
\newcommand{\LL}{{\rm \mathbf L}}
\newcommand{\kk}{{\rm \mathbf k}}
\newcommand{\pr}{p_R}
\newcommand{\pl}{p_L}
\newcommand{\qr}{q_R}
\newcommand{\ql}{q_L}
\newcommand{\ep}{\epsilon}
\newcommand{\om}{\omega}
\newcommand{\Om}{\Omega}
\newcommand{\Omb}{\bar{\Omega}}
\newcommand{\ga}{\gamma}

\newcommand{\GH}{\Gamma(H)}
\newcommand{\Si}{S^{-1}}
\newcommand{\Xbar}{\bar{X}}
\newcommand{\Ybar}{\bar{Y}}
\newcommand{\Zbar}{\bar{Z}}
\newcommand{\Lbar}{\bar{\Lambda}}

\newcommand{\re}{{\,\hbox{$\textstyle\triangleright$}\,}}
\newcommand{\li}{{\,\hbox{$\textstyle\triangleleft$}\,}}

\newcommand{\arr}{\!\rightharpoonup\!}
\newcommand{\arl}{\!\leftharpoonup\!}
\newcommand{\id}{{\rm id}}
\newcommand{\idM}{\id_M}
\newcommand{\one}{{\mathbf 1}}
\newcommand{\e}{{\mathbf 1}}

\newcommand{\bra}{\langle}
\newcommand{\ket}{\rangle}
\def\End{\mbox{End}\,}

\def\Alg{\mbox{Alg}}
\def\Hom{\mbox{Hom}\,}

\def\Aut{\mbox{Aut}\,}
\def\Mat{\mbox{Mat}\,}
\newcommand{\tp}{\otimes}
\newcommand{\tph}{\otimes_H}
\newcommand{\reli}{\bowtie}

\def\cros{\,\raise1.9pt\hbox{$\scriptscriptstyle  > $}\!
          \raise1.5pt\hbox{$\scriptstyle\triangleleft$}\,}
\def\>cros{\cros}
\def\<cros{\,\raise1.5pt\hbox{$\scriptstyle\triangleright$}\!
           \raise1.9pt\hbox{$\scriptscriptstyle < $}\,}

\newcommand{\res}{\upharpoonright }
\newcommand{\cop}{\Delta}
\newcommand{\bcop}{\bar{\Delta}}
\newcommand{\Eqref}[1]{Eq. \eqref{#1}}
\newcommand{\coh}{^{coH}}
\newcommand{\imn}{i_{MN}}

\newcommand{\0}{_{(0)}}
\newcommand{\1}{_{(1)}}
\newcommand{\2}{_{(2)}}
\newcommand{\3}{_{(3)}}

\begin{document}
\title{Integral Theory for Quasi-Hopf Algebras}
\author{Frank Hau{\ss}er}
\address{Frank Hau{\ss}er \\ Universit{\`a} di Roma 
"La Sapienza", 
        Dipartimento di Matematica\\
        P.le Aldo Moro 2, 00185 Roma, Italia}
\email{hausser@mat.uniroma1.it} 
\thanks{F.Hausser supported by EU TMR Network {\em Noncommutative Geometry}}

\author{Florian Nill}
\address{Florian Nill \\St{\"u}cklenstr. 1a\\ D-81247 M{\"u}nchen, Germany}
\email{nill@physik.fu-berlin.de}
\thanks{F.Nill supported by Arbeitsamt M\"unchen under KuNr. 843A416150}

\begin{abstract}
We generalize the fundamental structure Theorem on Hopf (bi)\-modules by
Larson and Sweedler to quasi-Hopf algebras $H$. For $\dim H<\infty$ this proves the existence and
uniqueness (up to scalar multiples) of integrals in $H$. Among other applications we prove a
Maschke-type Theorem for diagonal crossed products as constructed by
the authors in \cite{HN1,HN2}.
\end{abstract}

\maketitle

\tableofcontents 

\section{Introduction}
A (left) integral $l$ in a (quasi-) Hopf algebra $H$ is an
element of $H$ satisfying for all $a \in H$ ($\ep$ denoting the counit)
\begin{equation*}
  al = \ep(a) l.
\end{equation*}
For finite dimensional Hopf algebras $H$ with
dual $\dH$, one may identify $H \equiv \widehat{\hat H} $ and a left
integral may equivalently be characterized as a functional on $\dH$
being invariant under the canonical right action of $H$ on $\dH$, i.e.
satisfying $l(\psi \arl a) = \ep (a) \, l(\psi), \quad \forall  a \in
H, \, \psi \in \dH$. 
Choosing $H$ to be the group algebra of a finite
group $G$, $\dH$ becomes the algebra of functions on $G$ (under
pointwise multiplication) and a (left) integral in $H$ yields  
a (right) translation invariant measure on $G$ (in this case $l$ is given by
the sum over all elements of $G$ and is also a right integral in $H$).
Thus the notion of integrals in Hopf algebras is a natural
generalization of the notion of translation invariant group measures.
   
It is wellknown \cite{Abe, Sweedler}, that a (left) integral in a finite
dimensional Hopf algebra always exists and is unique up to
normalization. The proof in \cite{Abe, Sweedler, LS} uses the 
fundamental structure Theorem on Hopf modules by Larson and
Sweedler \cite{LS}.  
A (right) Hopf $H$-module  $M$ of a Hopf algebra $H$  
is a (right)  $H$-module which is also a (right) $H$-comodule such that the
corresponding $H$-action ($\cdot$) and $H$-coaction ($\rho$) are
compatible in a natural way (i.e. $\rho (m \cdot a) =  
\rho(m) \cdot \cop (a), \quad m \in M, \, a \in H$). 
Denoting the subspace of {\em coinvariants} of $M$ (i.e. of elements
satisfying $\rho (m) = m \tp \e$) by $M^{coH}$,  
the structure Theorem states that 
\begin{equation}
\label{0.1}
M \cong M^{coH} \tp H  \quad \quad \text{as Hopf $H$-modules},
\end{equation}
where the Hopf $H$-module structure on $M^{coH} \tp H$ is the
trivially extended Hopf $H$-module structure on $H$ (given by right
multiplication and coproduct) and the isomorphism $M^{coH} \tp H \to
M$ is given by $m \tp x \mapsto m \cdot x$.

For finite dimensional Hopf algebras the isomorphism \eqref{0.1} 
immediately implies the existence and uniqueness of (left)
integrals, since in this 
case also the dual $\dH$ becomes a  Hopf $H$-module by
dualizing the Hopf module structure of $H$. Then the left integrals in
$\dH$ are precisely the coinvariants and due to finite dimensionality
the subspace of coinvariants has to be one-dimensional.
This proves existence and uniqueness of integrals in 
$\dH$, and therefore also in $H$, by 
interchanging the role of $H$ and $\dH$. 

In this paper we generalize the structure Theorem to
quasi-Hopf algebras as introduced by Drinfeld \cite{Dr90}, and show
that many results of the integral theory of Hopf algebras carry over
to the quasi-Hopf case.
We will recall the precise definition of
quasi-Hopf algebras in Section~2.
The main difference compared to
ordinary Hopf algebras is the weakening of the coassociativity axiom,
i.e. the coproduct $\cop$ satisfy 
\begin{equation}
  \label{0.2}
 (\id\tp\cop)(\cop(a))\, \phi = \phi\,(\cop\tp\id)(\cop(a)), \quad
  a\in H
\end{equation}
with invertible $\phi \in H\tp H\tp H$, whereas in the Hopf case $\phi$
would be trivial. In particular this 
implies that the dual space $\dH$ in general is not an algebra.  
Equation~\eqref{0.2} already suggests that there is
an appropriate notion of quasi-Hopf $H$-bimodules (the guiding
example being $H$ itself), whereas a
meaningful notion of quasi-Hopf $H$-modules does not exist. 
However the structure Theorem for Hopf modules may as well be
formulated for Hopf $H$-bimodules and we will generalize this version
in Section~3 to quasi-Hopf algebras. We then provide for finite
dimenional $H$ a quasi-Hopf $H$-bimodule structure on the dual $\dH$. As
in \cite {LS} this implies 
\begin{equation*}
\dH \cong \L \tp H \quad \text{as quasi-Hopf $H$-bimodules,}
\end{equation*}
where $\L$ is a one-dimensional subspace of $\dH$, whose elements will
be called (left) {\em cointegrals} on $H$ (recall that $\dH$ is no algebra,
therefore neither does it make sense to speak of integrals in $\dH$
nor is it possible to simply interchange the role of $H$ and $\dH$ to
arrive at the existence and uniqueness of integrals in $H$). We then
show that every nonzero cointegral is nondegenerate (as a functional
on $H$), i.e. every finite dimensional quasi-Hopf algebra is
Frobenius, which implies existence and uniqueness of (left) integrals.
The existence part has already been proved before by F.~Panaite
\cite{PO} using different methods, see below.

In Section~5 we develop a notion of Fourier transformations and show
that these are given in
terms of cointegrals by the same formula as for ordinary Hopf
algebras. We then give a characterization of semisimple quasi-Hopf
algebras which essentially is due to  F.~Panaite \cite{P3}. Moreover
we show that $H$ is unimodular iff the modular automorphism of any
(left) cointegral is given by the square of the antipode $S$, or
equivalently, if every left integral is $S$-invariant. In Section~6 we 
generalize Radford's formula for the fourth
power of the antipode to quasi-Hopf algebras. In Section~7 we provide
a characterization of cointegrals in terms of {\em cocentral bilinear
forms} on $H$. This will allow as to give the suitable substitute
for the notion of cosemisimplicity needed to prove in Section~8
a Machke type Theorem for diagonal crossed products $A\reli \dH$ as
constructed in \cite{HN1} - an important example being the quantum
double $D(H) = H \reli \dH$, see \cite{HN2}.
In fact this has been our original motivation to develop an integral
theory for quasi-Hopf 
algebras. 

We remark that there is an alternative - and in fact very short -
approach to prove the existence 
and uniqueness of integrals in finite dimensional Hopf algebras, by
providing a direct formula for a nonzero integral, this way 
avoiding the use of the structure Theorem, see \cite{Dae97}.  This
approach has been successfully generalized to quasi-Hopf algebras by
\cite{PO}, proving the existence of integrals. Yet the
generalization of the uniqueness result along the lines of
\cite{Dae97} seems to be out of reach.

%%% Local Variables: 
%%% mode: latex
%%% TeX-master: "integral"
%%% End: 

\section{Quasi-Hopf Algebras}
We recall some definitions and properties of quasi-Hopf algebras as
introduced by Drinfeld \cite{Dr90} and fix our notations. Throughout $\kk$ will
be a field and all algebras and linear spaces will be over $\kk$.  All
algebra morphisms are supposed to be unital.
The dual of a linear space $V$ is denoted $\hat V:=\Hom\!_\kk(V,\kk)$.

A {\em quasi-bialgebra} $(H,\cop,\ep,\phi)$ is a unital $\kk$-algebra $H$ with
algebra morphisms $\cop:\, H\to H\tp H$ (the {\em coproduct}) and $\ep: \,
H\to \kk$ (the {\em counit}), and an invertible element $\phi \in H\tp H\tp H$
(the {\em reassociator}), such that
\begin{gather}
  \label{1.1}
  (\id\tp\cop)(\cop(a))\, \phi = \phi\,(\cop\tp\id)(\cop(a)), \quad
  a\in H, \\
   \label{1.2}
   (\id\tp\id\tp\cop)(\phi)\,(\cop\tp\id\tp\id)(\phi) =
   (\e\tp\phi)\,(\id\tp\cop\tp\id)(\phi)\,(\phi\tp\e), \\
   \label{1.3}
   (\ep\tp\id)\circ\cop = \id = (\id\tp\ep)\circ\cop,  \\
   \label{1.4}
   (\id\tp\ep\tp\id)(\phi)=\e\tp\e.
\end{gather}
The identities \eqref{1.1}-\eqref{1.4} also imply $(\ep\tp\id\tp\id)(\phi) =
(\id\tp\id\tp\ep)(\phi)= \e\tp\e$.  As for Hopf algebras we denote $\cop(a) =
\sum_i a^i_{(1)}\tp a^i_{(2)}= a\1\tp a\2 $, but since $\cop$ is only
quasi-coassociative we adopt the further convention
\begin{equation*}
  (\cop\tp\id)(\cop (a)) = a_{(1,1)}\tp a_{(1,2)} \tp a_{(2)},
   \quad (\id\tp\cop)(\cop(a)) = a_{(1)}\tp
  a_{(2,1)}\tp a_{(2,2)}.
\end{equation*}
Furthermore we denote $\cop^{op}(a) := a\2\tp a\1$ and 
\begin{equation*}
\phi=  X^j\tp Y^j \tp Z^j;\quad  \phi^{-1} = \Xbar^j\tp \Ybar^j\tp \Zbar^j,
\end{equation*}
suppressing the summation symbol $\sum_j$.

A quasi-bialgebra $H$ is called {\em quasi-Hopf algebra}, if there is a
linear antimorphism $S: \, H \rightarrow H$ (the {\em antipode}) and elements
$\alpha, \beta \in H$ satisfying for all $a\in H$
\begin{align}
  \label{1.7}
  & S(a_{(1)})\alpha a_{(2)} = \alpha \ep (a), \,\,\quad
  a_{(1)} \beta S(a_{(2)}) = \beta \ep(a),\\
  \label{1.8}
  & X^j \beta S(Y^j)\alpha Z^j = 1 = S(\Xbar^j)\alpha \Ybar^j \beta
  S(\Zbar^j).
\end{align}
We will always assume $S$ to be invertible. Equations \eqref{1.7},\eqref{1.8}
imply $\ep\circ S = \ep$ and $\ep (\alpha \beta) = 1$.
Thus, by rescaling $\alpha$ and $\beta$ we may without loss assume $\ep(\alpha)=\ep(\beta)=1$.

We will also frequently use the following notation: If $\psi = \sum_i \psi_i^1
\tp\dots \psi_i^m \,\in H^{\tp^m}$, then, for $m\le n$, $\psi^{n_1 n_2 \dots
  n_m} \in H^{\tp^n}$ denotes the element in $H^{\tp^n}$ having $\psi_i^k$ in
the ${n_k}^{\rm th}$ slot and $\e$ in the remaining ones.  

Next we recall that an invertible element $F\in H \tp H$  satisfying
$(\ep\tp\id)(F) = (\id\tp\ep)(F) = \e$, induces a so-called {\em twist
  transformation}
\begin{align}
  \notag
  \cop_F(a) : &= F \cop(a) F^{-1}, \\
   \label{1.6}
   \phi_F : &= (\e\tp F)\, (\id\tp\cop)(F) \, \phi \, (\cop\tp\id)(F^{-1})\,
   (F^{-1}\tp\e),
\end{align}
and $(H,\cop_F,\ep,\phi_F)$ is again a quasi-bialgebra.

It is well known that the antipode of a Hopf algebra is also an anti coalgebra
morphism, i.e. $\cop(a) = (S\tp S)\big(\cop^{op}(\Si(a))\big)$. For quasi-Hopf
algebras this is true only up to a twist.
 Following Drinfeld we define the
elements $\gamma, \delta \in H\tp H$ by 
\footnote{suppressing summation symbols}
\begin{align}
  \label{1.9'}
  \gamma & := (S(U^i)\tp S(T^i))\, (\alpha\tp\alpha)\, (V^i\tp W^i),
  \\
   \label{1.10}
   \delta & := (K^j\tp L^j)\, (\beta\tp\beta)\, (S(N^j)\tp S(M^j)),
\end{align}
where
\begin{align*}
  T^i\tp U^i \tp V^i \tp W^i &=(\e\tp \phi^{-1})\,(\id\tp\id\tp\cop)(\phi),
  \\  
  K^j \tp L^j \tp M^j \tp N^j & = (\phi\tp\e)\,
   (\cop\tp\id\tp\id)(\phi^{-1}). 
\end{align*}
With these definitions Drinfeld has shown in \cite{Dr90}, that $f\in H\tp H$
given by
\begin{equation}
  \label{1.11}
  f:= (S\tp S)(\cop^{op}(\Xbar^i)) \,\gamma \,
  \cop(\Ybar^i\beta S(\Zbar^i)) 
\end{equation}
is invertible and satisfies for all $a\in H$
\begin{equation}
  \label{1.12} f\cop(a)f^{-1} = (S\tp S)\big(\cop^{op}(\Si(a))\big).
\end{equation}
The elements $\gamma, \delta$ and the twist $f$ fulfill the relations
\begin{equation}
  \label{1.13}
  f\, \cop(\alpha) = \gamma,\quad \cop(\beta) \, f^{-1} = \delta\,.
\end{equation}
Furthermore, the corresponding twisted reassociator \eqref{1.6} is given by
\begin{equation}
  \label{1.17}
  \phi_f = (S\tp S\tp S)(\phi^{321}).
\end{equation}
Setting
\begin{equation}
  \label{1.14}
h : = (\Si\tp\Si)(f^{21}),
\end{equation}
the above relations imply
\begin{align}
  \label{1.15}
  h\cop(a) h^{-1} &= (\Si\tp\Si)\big(\cop^{op}(S(a))\big),  \\
   \label{1.18}
   \phi_h &= (\Si\tp\Si\tp\Si)(\phi^{321}).
\end{align}
Let $H_{op}$ denote the algebra with opposite multiplication, then
$H_{op}$ becomes a quasi-Hopf algebra by setting $\phi_{op}:=
\phi^{-1}$, $S_{op} := S^{-1}$, $\alpha_{op} := S^{-1}(\beta)$,
$\beta_{op} := S^{-1} (\alpha)$, implying $h_{op} = f^{-1}$.
 
We will also need a generalization of Hopf algebra formulae of the type $a\1
\tp a\2 S(a\3) = a \tp \e$ to the
quasi-coassociative setting. For this one uses the following elements
$\qr,\pr,\ql,\qr \in H\tp H$, see e.g.  \cite{HN1,HN2},
\begin{align}
\label{1.20}
  \qr &:= X^i\tp \Si(\alpha\,Z^i)\,Y^i ,
     &  \pr &: = \Xbar^i\tp \Ybar^i\,\beta\, S(\Zbar^i) ,  \\
\label{1.21}
  \ql &:=S(\Xbar^i)\,\alpha\, \Ybar^i \tp \Zbar^i , & \pl &:=
    Y^i\,\Si(X^i\,\beta)\tp Z^i .
\end{align}
They obey the relations (for all $a \in H$)
\begin{align}
 \label{1.22}  
 [\e\tp \Si(a\2)]\,\qr\, \cop(a\1) &= [a\tp\e]\, \qr , \\
 \label{1.22'}  
 \cop(a\1)\, \pr\, [\e\tp S(a\2)] &= \pr\, [a \tp \e], \\
\label{1.23}
[S(a\1)\tp\e]\, \ql \, \cop(a\2)& = [\e\tp a]\, \ql,\\
 \label{1.23'}  
 \cop(a\2)\,\pl\, [\Si(a\1)\tp \e]& = \pl\, [\e\tp a],
\end{align}
and (writing $q_R = q^1_R \tp q^2_R$, etc. suppressing the summation symbol
and indices)
\begin{align}  
\label{1.24}
\cop(\qr^1) \, \pr\, [\e\tp S(\qr^2)] &=\e\tp\e,
& [\e\tp\Si(\pr^2)]\, \qr \, \cop(\pr^1)& = \e\tp\e, \\
\label{1.25}
\cop (\ql^2) \, \pl \, [\Si (\ql^1) \tp \e] &=\e\tp\e, & [S(\pl^1) \tp \e]\,
\ql \, \cop(\pl^2) & = \e\tp\e.
\end{align}
Moreover, the following formula has been proved in \cite{HN1}
\begin{align}
  & (q_R \tp \e) \, (\cop\tp \id)(q_R)\,
  \phi^{-1}\notag\\
\label{1.26}
&\quad\quad\quad = [\e \tp S^{-1}(Z^i)\tp S^{-1}(Y^i)] \, [\e \tp h]\, (\id
\tp \cop)\big(q_R \cop(X^i)\big).
\end{align}

\bsn
We remark that the algebraic properties of a quasi-Hopf algebra $H$ may be
translated into corresponding properties of its representation category
$\makebox{Rep}\, H$. More precisely $\makebox{Rep}\, H$ is a rigid monoidal
category \cite{Dr90}. This implies that also the elements
$f,h,q_R,q_L,p_R,p_L$ and their properties as stated above may be nicely
described and understood in categorical terms, i.e. they define natural
transformations, see \cite{HN1}. One may also
use a grafical description to obtain and describe their properties, see
\cite{HN2}.  

\bsn
We denote $\Gamma(H):=\Alg(H,\kk)$ the group of 1-dimensional representations of $H$, where for 
$\mu,\nu\in\GH$ we put $\mu\nu:=(\mu\tp\nu)\circ\cop$.
Note that for 1-dimensional representations this product is indeed
strictly associative with unit $\hat\one\equiv\ep$.
Clearly, for any $\mu\in\GH$ \Eqref{1.8} implies $\mu(\alpha)\neq 0 \neq\mu(\beta)$.
Putting $\bar\mu:=\mu\circ S$ we conclude $\bar\mu\mu=\mu\bar\mu=\ep$ by \eqref{1.7}
and therefore $\bar\mu=\mu^{-1}$.
Thus $\GH$ is indeed a group.

\bsn
Finally we introduce $\dH$ as the dual space of $H$ with its natural
`multiplication` $\dH\tp\dH \to \dH$, which however is no longer
associative
\begin{equation*}
  \bra \varphi\psi\mid a\ket : = \bra \varphi\tp \psi\mid \cop(a)\ket, \quad
  \bra \hat{\e} \mid a \ket : = \ep(a),
\end{equation*}
where $\varphi,\psi \in \dH,\,a \in H$ and where 
$\bra \cdot \mid \cdot \ket : \dH \tp H \to \kk$ denotes
the dual pairing.
We have $\hat{\e} \varphi = \varphi \hat{\e} = \varphi$.
Transposing the right and left multiplication on $H$ one
obtains on  $\dH$ the left and right $H$-actions  
\begin{equation} \label{1.30}
  \bra a \arr \varphi\mid b\ket : =  \bra \varphi\mid ba \ket, \quad 
  \bra\varphi \arl a\mid b\ket :=\bra \varphi\mid ab\ket,
  \quad a,b \in H, \, \varphi \in \dH 
\end{equation}
satisfying
$
  a \arr (\varphi \psi) = (a\1 \arr \varphi)(a\2 \arr \psi)
$
and
$(\varphi \psi) \arl a = (\varphi \arl a\1)(\psi \arl a\2).
$
For $\varphi\in\dH$ and $a\in H$ we also use the dual notation
$$
\varphi\arr a:=a\1\bra\varphi\mid a\2\ket\quad,
\quad a\arl\varphi:=\bra\varphi\mid a\1\ket a\2\,.
$$
%For $\varphi\in\Gamma(H)\subset\dH$ this defines a
%$\Gamma(H)$-bimodule structure on $H$. 
If $H$ is finite dimensional, then $\dH$ is also equipped with a
coassociative coalgebra structure $(\hat{\cop},\hat{\ep})$ given by $\bra
\hat{\cop}(\varphi) \mid a \tp b\ket:=\bra \varphi \mid ab\ket$ and
$\hat{\ep}(\varphi) : = \bra \varphi \mid \e \ket$.
 In this case we have
\begin{equation*}
  a \arr \varphi = \varphi_{(1)} \bra \varphi_{(2)}\mid a \ket, \quad 
  \varphi \arl a = \varphi_{(2)} \bra \varphi_{(1)}\mid a\ket,
  \quad a \in H, \, \varphi \in \dH 
\end{equation*}
 and $\hat{\cop}(\varphi\psi) = \hat{\cop}(\varphi) \hat{\cop}(\psi)$.

%%% Local Variables: 
%%% mode: latex
%%% TeX-master: "integral"
%%% End: 

\section{Quasi-Hopf Bimodules} 
Throughout let $H$ denote a quasi-Hopf algebra
over the field $\kk$. As a natural generalization of the notion of Hopf
bimodules we define
\begin{definition}
\label{def2.1}
Let $M$ be an $H$-bimodule and let $\rho :M\rightarrow M\otimes H$ be an
$H$-bimodule map.  Then $(M,\rho )$ is called a {\em quasi-Hopf $H$-bimodule},
if
\begin{align}
  (\idM\otimes \ep)\circ \rho  &= \idM , \label{2.1} \\
  \phi \cdot (\rho \otimes \idM)(\rho (m)) &= (\idM\otimes \Delta )(\rho
  (m))\cdot \phi ,\quad \forall m\in M.
\label{2.2}
\end{align}
\end{definition}        
We call $\rho $ a {\em right $H$-coaction} on $M$ and property \eqref{2.2} the 
{\em quasi-coassociativity} of $\rho$. 
Similarly as for coproducts we use the suggestive notation 
\begin{equation*}
  \rho(m) = m\0 \tp m\1, \quad (\rho\tp \idM)(\rho(m)) = m_{(0,0)}\tp
  m_{(0,1)} \tp m\1, \quad \text{etc.}
\end{equation*}
More specifically, one may also call $(M,\rho)$ a {\em right} quasi-Hopf
$H$-bimodule, whereas a {\em left} quasi-Hopf $H$-bimodule $(M,\Lambda)$ would
be defined analogously, except that $\Lambda:M\to H\tp M$ would now be a
quasi-coassociatve left coaction. If not mentioned explicitely, we will always
work with right coactions.

A trivial example is of course given by $M=H$ and $\rho =\Delta $. More
generally, one straightforwardly checks the following
\begin{lemma}
\label{lem2.2}
Let $(V,\re)$ be a left $H$-module. Then $V\tp H$ becomes a quasi-Hopf
$H$-bimodule by putting for $a,b,x\in H$ and $v\in V$
\begin{align*}
  a\cdot (v\otimes x)\cdot b & := (a_{(1)}\re v)\otimes a_{(2)}xb, \\
  \rho _{V\otimes H}(v\otimes x) & := \bar X^{i}\re v\otimes \bar
  Y^{i}x_{(1)}\otimes \bar Z^{i}x_{(2)}.
\end{align*}
\end{lemma}
In the special case $V=\kk$ with $H$-action given by some $\gamma\in\GH$ we may identify
$V\tp H\cong H$ to obtain a quasi-Hopf $H$-bimodule structure on $H$, 
denoted by $(H_\gamma,\rho_\gamma)$.
This leads to
\begin{corollary}
\label{cor2.3}
Let $\gamma\in\GH$ and put $T_\ga:=\ga(\Xbar^i)\Ybar^i\tp\Zbar^i\in H\tp H$.
Denote $H_\ga:=H$ considered as an $H$-bimodule with actions
\begin{equation*}
a\cdot x\cdot b :=(a\arl\ga)xb,\quad a,b\in H,\,x\in H_\ga .
\end{equation*}
Put $\rho_\ga:H_\ga\to H_\ga\tp H,\ \rho_\ga(x):=T_\ga\cop(x)$.
Then $(H_\ga,\rho_\ga)$ provides a quasi-Hopf $H$-bimodule.
\end{corollary}
Our aim is to generalize the fundamental structure Theorem on Hopf
(bi-)\-modules 
by Larson and Sweedler \cite{LS} to quasi-Hopf algebras $H$.  Thus, any
quasi-Hopf $H$-bimodule $M$ will be shown to be isomorphic to some $V\tp H$ as
in Lemma~\ref{lem2.2}, where $V\equiv M^{coH}\subset M$ will be a suitably
defined subspace of coinvariants.

To this end we first provide what will be a projection $E:M\to M^{coH}$, which
for ordinary Hopf algebras $H$ would be given by $E(m)=m\0 \cdot S(m\1)$.  In
this case $E(M)\equiv M^{coH}$ would be invariant under the adjoint $H$-action
$a\re m:=a\1\cdot m\cdot S(a\2)$ and we would have $E(a\cdot m)=a\re E(m),\ 
\forall a\in H,\,m\in M$.

In the quasi-Hopf case we first appropriately generalize this last property.
Here the basic idea is that the "would-be-adjoint" action $\re$ should satisfy
$ a\cdot m = (a\1\re m)\cdot a\2$. In fact, we will see that this only holds
for $m\in M^{coH}$. For $m\in M$ and $a\in H$ we now define
\begin{align}
  E(m) &:= \qr^1 \cdot m\0 \cdot \beta S(\qr^2 m\1)
  \label{2.3}\\
  a\re m &:= E(a\cdot m) \label{2.4},
\end{align}
where $\qr^1\tp\qr^2 \equiv \qr \in H\tp H$ is defined in \eqref{1.20}.
\begin{proposition}
\label{prop2.3}
Let $M$ be a quasi-Hopf $H$-bimodule and let $E$ and $\re$ be given as above.
Then for all $a,b\in H$ and $m\in M$
\begin{equation*}
\begin{aligned}
 \text{\rm (i)} & & E(m\cdot a) &= E(m)\,\ep(a)       \\
 \text{\rm (ii)} & &E^2 &= E                         \\
  \text{\rm (iii)} && a\re E(m) &= E(a\cdot m) \equiv a\re m       \\
  \text{\rm (iv)} && (ab)\re m &=  a\re (b\re m)       \\
  \text{\rm (v)} &&  a\cdot E(m) &=[a\1\re E(m)]\cdot a\2  \\
  \text{\rm (vi)}  & &\quad E(m\0)\cdot m\1 &= m.\\
  \text{\rm (vii)}  & &\quad E(E(m)\0)\tp E(m)\1 &= E(m)\tp\one.
\end{aligned}
\end{equation*}
\end{proposition}
\begin{proof}
  All properties (i) - (vii) follow easily from the properties of quasi-Hopf
  algebras as stated in Section~2. We will give a rather detailed proof
  such that the unexperienced reader may get used to the techniques used when
  handeling formulae involving iterated non coassociative coproducts.
  We will in the following denote $q_R = q = q^1\tp q^2$.
  Equality (i) follows directly from the antipode property $a\1\beta S(a\2) =
  \ep(a)\beta$. To show (ii) one uses (i) to compute
  \begin{align*} 
    E^2 (m) & = E\big(q^1 \cdot m\0 \cdot\beta S(q^2 m\1)\big) \\
      & = E(q^1 \cdot m\0) \, \ep\big(\beta S(q^2 m\1)\big) \\
     & = E(m) \, \ep(\alpha \beta) = E(m)
  \end{align*}
by \eqref{1.4} and \eqref{2.1}.
Equality (iii) is obtained similarly: 
\begin{equation*}
  a\re E(m) = E\big(a\cdot E(m)\big) 
            =E\big(aq^1 \cdot m\0 \cdot \beta S(q^2 m\1)\big)
            =E(a \cdot m).
\end{equation*}
Property (iv) follows immediately from (iii).
To show (v) we use \eqref{1.22}:
\begin{align*}
  a \cdot E(m) &= aq^1 \cdot m\0 \cdot \beta S(q^2 m\1) \\
  &= q^1 a_{(1,1)} \cdot m\0 \cdot \beta S(q^2 a_{(1,2)} m\1 )\, a\2\\
   &= E(a\1 \cdot m)\cdot a\2 \\
   &= \big[a\1 \re E(m)\big] \cdot a\2 .
\end{align*}
Part (vi) follows by using formula \eqref{1.20} for $q = q_R$, \eqref{2.2} and
then \eqref{1.7},\eqref{1.8}
\begin{align*}
  E(m\0)\cdot m\1 & = X^i \cdot m_{(0,0)} \cdot \beta S(Y^i m_{(0,1)})
  \alpha Z^i m\1 \\
  &= m \cdot X^i \beta S(Y^i) \alpha Z^i = m.
\end{align*}
Finally, we prove (vii) by using part (i), Eq. \eqref{2.2} and the left identity in 
\eqref{1.24} to compute
\begin{align*}
E(E(m)\0)\tp E(m)\1 &=E(q^1\1\cdot m_{(0,0)})\tp q^1\2m_{(0,1)}\beta S(q^2 m\1)\\
&=E(q^1\1\Xbar^i\cdot m_{(0)})\tp q^1\2\Ybar^i m_{(1,1)}\beta 
    S(q^2 \Zbar^i m_{(1,2)})\\
&=E(q^1\1 \Xbar^i\cdot m)\tp q^1\2 \Ybar^i \beta S(\Zbar^i )S(q^2)\\
&=E(m)\tp\one\,.
\end{align*}
\end{proof}
Due to part (ii), (vi) and (vii) of Propostion~\ref{prop2.3} the following notions 
of {\em coinvariants} all coincide
\begin{definition}
\label{def2.4}
The space of {\em coinvariants} of a quasi-Hopf $H$-bimodule $M$ is defined to
be 
\begin{equation*}
  \McoH:=E(M)\equiv\{n\in M\mid E(n)=n\}\equiv\{n\in M\mid E(n\0)\tp
  n\1=E(n)\tp\one\} 
\end{equation*}
\end{definition}
Let us first check this definition for the type of quasi-Hopf $H$-bimodules
described in Lemma~\ref{lem2.2}.
\begin{lemma}
\label{lem2.5}
The coinvariants of the quasi-Hopf $H$-bimodule $V\tp H$ in Lemma \ref{lem2.2}
are given by
$
(V\tp H)^{coH} =V\tp\one,
$
and for $v\in V,\ x\in H,$ we have
$
E(v\tp x)=v\tp\ep(x)\one.
$
\end{lemma}
\begin{proof}
  The definitions in Lemma~\ref{lem2.2} imply $(v\tp x) = (v\tp \e)
  \cdot x$, which by Proposition~\ref{prop2.3} (i) further implies
  \begin{equation*}
    E(v\tp x) = E(v\tp \e) \, \ep(x).
 \end{equation*}
 Thus we are left to show the identity
    $E(v\tp \e) = v \tp \e$ (denoting $q = q_R$):
 \begin{align*}
      E(v\tp \e) &= q^1 \cdot \big[ \Xbar^i \re v \tp \Ybar^i\big] \cdot
      \beta S(q^2 \Zbar^i) \\
     & = q^1\1 \Xbar^i \re v \tp q^1\2 \Ybar^i \beta S(\Zbar^i)
      S(q^2) \\
     & = v\tp \e ,
  \end{align*}
where the last equality follows from \eqref{1.24}.
\end{proof}
For later purposes let us also consider the $\kk$-dual $\hat M$ as a left
$H$-module by transposing the right $H$-action on $M$. Denote $\hat M^H$ the
invariants under this action, i.e.
\begin{equation*}
\hat M^H:=\{\psi\in\hat M\mid \psi(m\cdot a)=\psi(m)\ep(a),\ \forall m\in
M,\,a\in H\} .
\end{equation*}
\begin{lemma}
\label{lem2.6} 
Let $E^T:\hat M\to\hat M$ denote the transpose projection of $E$.  Then
\begin{equation*}
E^T(\hat M)=\hat M^H.
\end{equation*}
\end{lemma}
\begin{proof}
  Clearly, part (i) of Proposition~\ref{prop2.3} implies $E^T(\hat
  M)\subset\hat M^H$. Conversely, if $\psi\in\hat M^H$ then \eqref{2.3}
  implies $\psi\circ E=\psi$, hence $\psi\in E^T(\hat M)$.
\end{proof}
We are now in the position to provide the fundamental Theorem for quasi-Hopf
$H$-bimodules.
By a morphism $(M,\rho)\to(M',\rho')$ of quasi-Hopf $H$-bimodules we mean
 an $H$-bimodule map $f:M\to M'$ satisfying
$\rho'\circ f=(f\tp\id)\circ\rho$.
\begin{theorem}
\label{thm2.7}  
Let $M$ be a quasi-Hopf $H$-bimodule. Consider $N\equiv\McoH$ as a left
$H$-module with $H$-action $\re$ as in \eqref{2.4}, and $N\tp H$ as a
quasi-Hopf $H$-bimodule as in Lemma~\ref{lem2.2}.  Then
\begin{equation*}
\nu: N\tp H\ni n\tp a\mapsto n\cdot a\in M
\end{equation*}
provides an isomorphism of quasi-Hopf $H$-bimodules with inverse given by
\begin{equation*}
\nu^{-1}(m) = E(m\0)\tp m\1\,.
\end{equation*}
\end{theorem}  
\begin{proof}
Using Proposition~\ref{prop2.3} (i), (vii) and $n=E(n)$ we compute
 \begin{align*}
   \nu^{-1} \circ \nu\, (n \tp a) & = E(n\0\cdot a\1) \tp n\1 a\2 \\
                   & = E(n\0) \tp n\1 \, a \\
                   & = E(n) \tp a = n \tp a. 
 \end{align*}
 Conversely $ \nu \circ \nu^{-1} (m) = E(m\0)\cdot m\1 = m$
 by Proposition~\ref{prop2.3} (vi). Thus $\nu$ is indeed an isomorphism of
 vector spaces. 
 
 We are left to show that $\nu$ also respects the quasi-Hopf $H$-bimodule
 structures.  By definition we have $a \cdot (n\tp x )\cdot b = E(a\1 \cdot n )
 \tp a\2 x b$ and therefore
\begin{align*}
  \nu\, \big(a \cdot (n\tp x) \cdot b \big) &= E(a\1 \cdot n) \cdot a \2 x b \\
        &= \big[ a\1 \re E(n)\big] \cdot a\2 x b\\
       & = a \cdot E(n) \cdot xb = a \cdot n\cdot x b\\
       & = a \cdot \nu (n \tp x ) \cdot b. 
\end{align*}
Here we have used Proposition~\ref{prop2.3} (iii) and (v) in
the second and third line, respectively. Thus $\nu$ is an $H$-bimodule map.
Finally, we show that $\nu^{-1}$ (and therefore $\nu$) are also $H$-comodule maps.
 \begin{align*}
 \rho_{N\tp H}\,(\nu^{-1}(m)) &=E(\Xbar^i\cdot m\0)\tp\Ybar^i m_{(1,1)}\tp\Zbar^i m_{(1,2)}\\
&=E(m_{(0,0)})\tp m_{(0,1)}\tp m\1\\
&=(\nu^{-1}\tp\id)(\rho_M(m))\,.
 \end{align*}
Here we have used part (iii) of Proposition~\ref{prop2.3} in the first line and 
part (i) together with the quasi-coassociativity \eqref{2.2} in the second line.
\end{proof}
\begin{corollary}
\label{cor2.6} 
For any quasi-Hopf $H$-bimodule $M$ we have
\begin{equation*} 
 M\coh = \{n\in M\mid\rho(n)=(\bar X^i\re n)\cdot\bar Y^i \tp \bar Z^i\}
\end{equation*}
\end{corollary}
\begin{proof}
If $\rho(n)=(\bar X^i\re n)\cdot\bar Y^i \tp \bar Z^i$ then by Proposition~\ref{prop2.3}(i) and (vi)
$E(n)=E(n\0)\cdot n\1=n$, whence $n\in N:=M\coh$. The inverse implication follows from
 $\rho_M=(\nu\tp\id)\circ\rho_{N\tp H}\circ\nu^{-1}$ and $\nu^{-1}(N)=N\tp\one$.
\end{proof}
Theorem~\ref{thm2.7} shows that there is a one-to-one correspondence (up to
equivalence) between left $H$-modules and quasi-Hopf $H$-bimodules.
 As for
ordinary Hopf algebras, this is actually an equivalence of monoidal
categories. To see this we note that if $M$ and $N$ are quasi-Hopf
$H$-bimodules, then so is $M\tp_H N$ with its natural $H$-bimodule structure, 
the coaction being given by
\begin{equation*}
\rho_{M\tp_H N}(m\tp n):=(m\0\tp n\0)\tp m\1 n\1\,.
\end{equation*}
In this way the category ${}_H\M_H^H$ of quasi-Hopf $H$-bi\-modules becomes a
strict mono\-idal category with unit object given by $H$. Moreover, we have
\begin{lemma}
\label{lem2.8}
Let $M$ and $N$ be quasi-Hopf $H$-bimodules. Then the map
\begin{equation*}
i_{MN}:\McoH\tp N^{coH}\ni m\tp n\mapsto (X^i\re m)\tp_H(Y^i\re n)\cdot Z^i \in
(M\tp_H N)^{coH}
\end{equation*}
provides an isomorphism of left $H$-modules.
\end{lemma}
\begin{proof}
Denote $E_M\,,E_N$ and $E_{M\tph N}$ the projections onto the
corresponding coinvariants.
To prove that $i_{MN}:M\coh\tp N\coh\to(M\tph N)\coh $ is bijective we claim
\begin{align}
\label{a}
\bar \imn\circ\imn &= \id \\
\label{b}
\imn\circ\bar\imn  &= E_{M\tph N},
\end{align}
where $\bar\imn:M\tph N\to M\coh\tp N\coh$ is  given by
\begin{equation}
\label{c}
\bar\imn(m\tph n)  := E_M(m\0)\tp E_N(m\1\cdot n)\,.
\end{equation}
By Proposition~\ref{prop2.3}(i) $\bar\imn$ is well defined and we have
\begin{equation}
\label{d}
\bar\imn(m\tph n) = m\tp E_N(n),\quad\forall m\in M\coh\,
\end{equation}
thus proving \eqref{a}.
To prove \eqref{b} let $m\in M\coh$ and $n\in N$. Then ($q:=q_R$)
\begin{align*}
E_{M\tph N}(m\tph n) &=
q^1\cdot m\0\tph n\0\cdot\beta S(q^2 m\1 n\1)
\\
& = q^1\cdot(\Xbar^i\re m)\cdot\Ybar^i\tph n\0\cdot\beta S(q^2\Zbar^i n\1)
\\
&= q^1\1 \Xbar^i\re m\tph q^1\2 \Ybar^i\cdot n\0\cdot\beta S(q^2\Zbar^i n\1)
\\
&= X^k\re m\tph q^1 Y^k\1\cdot n\0\cdot\beta S(q^2 Y^k\2n\1)Z^k
\\
&= X^k\re m\tph (Y^k\re n)\cdot Z^k\,,
\end{align*}
where in the second line we have used Corollary~\ref{cor2.6}, in the third line
 Proposition~\ref{prop2.3}(v) and in the fourth line
the identity
\begin{equation*}
(\cop\tp\id)(q_R)\phi^{-1} = [X^k\tp\one\tp S^{-1}(Z^k)][\one\tp q_R\cop(Y^k)]
\end{equation*}
which follows easily from \eqref{1.2} and \eqref{1.7}. 
Thus, by Proposition~\ref{prop2.3}(vi) and (iii) we conclude for general $m\in M, n\in N$
\begin{align*}
E_{M\tph N}(m\tph n) &=
E_{M\tph N}(E(m\0)\tph m\1\cdot n) 
\\
&= X^k\re E(m\0)\tph [Y^k\re (m\1\cdot n)]\cdot Z^k
\\
&= (\imn\circ\bar\imn)(m\tph n)\,.
\end{align*}
This proves that $\imn$ is bijective.
To prove that it is $H$-linear we compute for $a\in H,\ m\in M\coh$ and $n\in N\coh$
\begin{align*}
a\cdot\imn(m\tp n) &= a\cdot(X^i\re m)\tph(Y^i\re n)\cdot Z^i
\\
&= a\1X^i\re m\tph(a_{(2,1)}Y^i\re n)\cdot a_{(2,2)}Z^i
\\
&= \imn(a_{(1,1)}\re m\tp a_{(1,2)}\re n)\cdot a\2\,.
\end{align*}
$H$-linearity of $\imn$ follows by applying $E_{M\tp_H N}$ to both sides
 and using Proposition~\ref{prop2.3} (i).
\end{proof}
Denoting the category of left $H$-modules by ${}_H\M$ Lemma~\ref{lem2.8} leads
to
\begin{proposition}
\label{prop2.9} 
There is an equivalence of monoidal categories ${}_H\M_H^H\cong{}_H\M$ given
on the objects by $M\mapsto \McoH$ and on the morphisms by $f\mapsto
f\res\McoH$.
\end{proposition}
\begin{proof}
If $f:M\to N$ is a morphism of quasi-Hopf $H$-bimodules, then $f(M\coh)\subset N\coh$,
$f\res M\coh$ is $H$-linear and
\begin{equation*}
f(m)=f(E(m\0)\cdot m\1)=f(E(m\0))\cdot m\1\,.
\end{equation*}
Thus, $f$ is uniquely determined by its restriction $f\res M\coh$, and by Theorem~\ref{thm2.7}
$(M\mapsto M\coh\,,\,f\mapsto f\res M\coh)$ provides an equivalence of categories
with reverse functor given by $V\mapsto V\tp H$ and $(f:V\to W)\mapsto (f\tp\id: V\tp H\to W\tp H)$.
By Lemma~\ref{lem2.8} these functors preserve the monoidal structures provided
we show that for any three objects $M,N,K$ in $\HMH$ the following
diagram commutes
\begin{equation}  
\label{diagr}
  \unitlength0.5cm
 \begin{picture}(24,7)
  \put(12,1){\makebox(0,0){$(M\tph N\tph K)\coh$}}
  \put(4,2.25){\vector(3,-1){4}}
  \put(4.4,1){\makebox(0,0){$i_{(M\tph N)K}$}}
  \put(19.5,2.25){\vector(-3,-1){4}}
  \put(20,1){\makebox(0,0){$i_{M(N\tph K)}$}}
  \put(4,3){\makebox(0,0){$(M\tph N)\coh\tp K\coh$}}
  \put(20,3){\makebox(0,0){$M\coh\tp(N\tph K)\coh$}}
  \put(4,5.5){\vector(0,-1){2}}
  \put(2.25,4.5){\makebox(0,0){$\imn\tp\id$}}
  \put(20,5.5){\vector(0,-1){2}}
  \put(21.75,4.5){\makebox(0,0){$\id\tp i_{NK}$}} 
  \put(4,6){\makebox(0,0){$(M\coh\tp N\coh)\tp K\coh$}}
  \put(9,6){\vector(1,0){5}}
  \put(11.5,6.5){\makebox(0,0){$\phi$}}
  \put(20,6){\makebox(0,0){$M\coh\tp(N\coh\tp K\coh)$}}
 \end{picture}
\end{equation}
To this end let $m\in M\coh,\,n\in N\coh,\,k\in K\coh$, then
\begin{align*}
\big(i_{(M\tph N)K} &\circ(\imn\tp\id)\big)(m\tp n\tp k) 
\\ 
 & = \imn(X^i\1\re m\tp X^i\2\re n)\tph(Y^i\re k)\cdot Z^i
\\
& = X^jX^i\1\re m\tph Y^jX^i\2\re n\tph (Z^j\1Y^i\re k)\cdot Z^j\2Z^i
\\
& = X^jX^i\re m\tph X^k Y^j\1Y^i\re n\tph (Y^kY^j\2Z^i\re k)\cdot Z^kZ^j
\\
& = X^jX^i\re m\tph i_{NK}(Y^j\1Y^i\re n\tp Y^j\2Z^i\re k)\cdot Z^j
\\
& = \big(i_{M(N\tph K)}\circ(\id\tp i_{NK})\big)(X^i\re m\tp Y^i\re
 n\tp Z^i\re k)\,, 
\end{align*}
where we have used Lemma~\ref{lem2.8} in the first line, $H$-linearity of $\imn$
in the second line, \eqref{1.2} in the third line and again
$H$-linearity of $i_{NK}$ in the last line.
\end{proof}
Note that $\HMH$ is strictly monoidal whereas $\HM$ is not.
It will be shown elsewhere that 
${}_H\M_H^H$ naturally coincides with the category of representations of the
two-sided crossed product $A:=H\>cros\hat H\<cros H$  constructed in \cite{HN1}.  Thus, by
methods of P.~Schauenburg \cite{Sch}, $A$ becomes a
$\times_H$-bialgebra (also called quantum
groupoid) in the sense of Takeuchi \cite{T}.
If $H$ is finite dimensional and Frobenius-separable (i.e. admits a
nondegenerate functional of index one), then by \cite{NSch} the Takeuchi quantum groupoids
are the same as the weak Hopf algebras of
\cite{BSz, N, BNS}.  This will approach a proof of an
anouncement of \cite{N}, saying that \footnote{over an algebraically closed field
  $\kk$ of characteristic zero} to any f.d. semi-simple quasi-Hopf algebra $H$
there is a (strictly coassociative!) weak Hopf algebra structure on $A:=H\>cros\hat H\<cros H$,
whose representation category obeys the same fusion rules.

\bsn
With the application to integral theory in mind we finally show, that for
finite dimensions any {\em left} quasi-Hopf $H$-bimodule
$(M,\Lambda)$ naturally gives rise to a dual {\em right} quasi-Hopf
$H$-bimodule $(M^*,\rho)$. 
As a linear space we put  $M^*=\hat M\equiv\Hom_\kk(M,\kk)$ with
$H$-bimodule structure given for $a,b\in H,\,m\in M$ and
$\psi\in\hat M$ by
\begin{equation}
\label{2.15}
\bra a\cdot \psi\cdot b\mid m\ket:=\bra\psi\mid S^{-1}(a)\cdot m \cdot
S(b)\ket\,. 
\end{equation}
To define the right $H$-coaction on $M^*$ we first deform the left
coaction $\Lambda$ on $M$ by putting $\Lbar: M\to H\tp M$ according to
\begin{equation}
\label{2.16}
\Lbar(m):= V\cdot\Lambda(m)\cdot U,
\end{equation}
where $U,V\in H\tp H$ are given by
\begin{align}
  U :=& f^{-1} (S\tp S)(q_R^{21}),
\label{2.17}\\
V :=&(S^{-1}\tp S^{-1})(p_R^{21})\,h,
\label{2.18}
\end{align}
the elements $f,h,q_R,p_R\in H\tp H$ being defined in
\eqref{1.11}-\eqref{1.21}.  With these definitions 
the linear map $\Lbar$ satisfies $(\ep\tp \idM) \circ \Lbar = \idM$ and 
\begin{gather}
\label{2.19}
  [\one\tp S^{-1}(a)]\cdot \Lbar(m) \cdot [\one\tp S(b)] = [a\2\tp\one]\cdot
  \Lbar\big(S^{-1}(a\1)\cdot m \cdot S(b\1)\big)\cdot[b\2\tp\one] \\
\label{2.20}
  [Y^i\tp Z^i\tp \one] \cdot (\id\tp\Lbar)\Lbar\big(S^{-1}(X^i)\cdot m\big)
  = (\cop\tp\id)\Lbar\big(m\cdot S(X^i)\big) \cdot [Y^i\tp Z^i\tp\one], 
\end{gather}
see Lemma~\ref{lem2.11} below. These identities imply the following
\begin{proposition}
\label{prop2.12} Let $(M,\Lambda)$
be a left quasi-Hopf $H$-bimodule and assume $H$ or $M$ finite dimensional. 
Consider $M^*$ as an $H$-bimodule as in
\eqref{2.15} and define $\rho :M^*\to M^*\tp H$ by identifying $M^*\tp
H\cong \Hom_\kk(M,H)$ and putting for $\psi\in M^*$ and $m\in M$
\begin{equation}\label{2.21}
[\rho(\psi)](m):=(\id_H\tp\psi)(\Lbar(m))\in H\,.
\end{equation}
Then $(M^*,\rho)$ provides a right quasi-Hopf $H$-bimodule.
\end{proposition}
\begin{proof}
  $H$-linearity of $\rho$ follows straightforwardly from \eqref{2.19}.
  Identifying 
\\ 
$\Hom\!_\kk(M, H\tp H)\cong M^*\tp H\tp H$,
quasi-coassociativity of $\rho$ is equivalent to the identity 
\begin{equation*}
  (Y^i\tp Z^i) \cdot \big[(\rho\tp\id)\rho(\psi)\big]\big(S^{-1}(X^i)\cdot
  m\big) 
  = \big[(\id\tp\cop)\rho(\psi)\big]\big(m\cdot S(X^i)\big) \cdot (Y^i\tp
  Z^i),
\end{equation*}
which follows from \eqref{2.20}.
\end{proof}
We conclude this section with the proof of the identities \eqref{2.19} and \eqref{2.20}.
\begin{lemma}
\label{lem2.11} 
For all $a,b\in H$ and $m\in M$ we have
  \begin{align}
  \label{2.22}
  U\,[\one\tp S(a)] &=  \Delta\big(S(a\1)\big)\, U \,[a\2\tp\one] \\
\label{2.23} [\one\tp
S^{-1}(a)]\,V &=  [a\2\tp\one]\,V\,\Delta\big(S^{-1}(a\1)\big) \\
\label{2.24}
\phi^{-1}\,(\id\tp\Delta)(U)\,(\one\tp U) &=
(\Delta\tp\id)\big(\Delta (S(X^i))U\big)\,(Y^i\tp Z^i\tp \one) \\
\label{2.25} 
(\Delta\tp\id)(V)\,\phi^{-1} = (Y^i\tp &Z^i\tp\one)(\one\tp
V)\, (\id\tp\Delta)\big(V\Delta\big(S^{-1}(X^i)\big)\big),
\end{align} 
implying \eqref{2.19} and \eqref{2.20}.
\end{lemma}
\begin{proof}
  Equation~\eqref{2.25} is equivalent to
  \begin{equation*}
  (\one\tp V)(\id\tp\Delta)(V)\phi = (\bar Y^i\tp \bar
  Z^i\tp\one)\,(\Delta\tp\id)\big(V\Delta (S^{-1}(\bar X^i) )\big).
  \end{equation*}
  Thus, noting that in $H_{op}$ the roles of $U$ and $V$ interchange,
  \eqref{2.23} and \eqref{2.25} reduce to \eqref{2.22} and
  \eqref{2.24}, respectively, in $H_{op}$.  
To prove  \eqref{2.22} we compute, using  \eqref{1.12} and \eqref{1.22} and
  denoting $q := q_R$
\begin{align*}
  \Delta\big(S(a\1)\big)\, U \,[a\2\tp\one] &= f^{-1}\,(S\tp S)\big(
  [S^{-1}(a\2) \tp
  \one ]\, q^{21} \, \cop^{op}(a\1)\big)\\
  & = f^{-1} \, (S\tp S)\big( [\one \tp a] \, q^{21}\big) \\
  & = U\,[\one\tp S(a)].
\end{align*}
To prove \eqref{2.24} we compute 
\begin{align*}
  \phi^{-1}&(\id\tp\Delta)(U)(\one\tp U)\\
  & \overset{\eqref{1.12}}{=}
  \phi^{-1}\,(\id\tp \cop)(f^{-1}) \, (\e\tp f^{-1}) \,
    (S\tp S\tp S)\big( (\e\tp q^{21})(\id\tp\cop^{op})(q^{21})\big) \\
  &\overset{\eqref{1.17}}{=} (\cop\tp\id)(f^{-1}) \, (f^{-1}\tp \e)\,
      (S\tp S\tp S)\big(\big[ (q\tp \e)(\cop\tp\id)(q)
  \phi^{-1}\big]^{321}\big) \\
   &\overset{\eqref{1.26}}{=}  (\cop\tp\id)(f^{-1}) \, (f^{-1}\tp \e)\,
  (S\tp S\tp S)\Big((\cop^{op}\tp\id)\big(q^{21}\,\cop^{op}(X^i)\big)\Big) \\
   &\qquad\qquad (f\tp\e)\, (Y^i\tp Z^i\tp \e)\\
  & \overset{\eqref{1.12}}{=}
    (\cop\tp \id)\big(\cop(S(X^i))\big) \, (\cop\tp\id)(f^{-1})  \,
     (\cop\tp\id)\big((S\tp S)(q^{21})\big)\\
   &\qquad\qquad (f^{-1}\tp\e)(f\tp\e) \,
     (Y^i\tp Z^i\tp \e)\\
  &\quad = (\cop\tp \id)\big(\cop(S(X^i))\big) \, (\cop\tp\id)(U)\, 
      (Y^i\tp Z^i\tp \e).
\end{align*}
Now \eqref{2.19} follows immediately from \eqref{2.22},\eqref{2.23} and
\eqref{2.20} from \eqref{2.24}, \eqref{2.25}.
\end{proof}

%%% Local Variables: 
%%% mode: latex
%%% TeX-master: "integral"
%%% End: 

\section{Integral Theory}
As in the original work of Larson and Sweedler \cite{LS}, the first
application of  
our previous results provides a theory of integrals and Fourier
transformations for finite dimensional quasi-Hopf algebras.
\begin{definition}
\label{def3.1}
An element $l\in H$ ($r\in H$) is called a {\em left (right) integral}, if
$al=\ep(a)l\, (ra=r\ep(a)),\ \forall a\in H$.  If $l$ is a left and a right
integral, then it is called {\em two-sided}. A left (right) integral $l$ is
called {\em normalized}, if $\ep(l)=1$. A {\em Haar integral} $e$ is a
normalized two-sided integral.
\end{definition}
Note that being the unit in the ideal of two-sided integrals a Haar integral
$e\in H$ is unique, provided it exists. In particular $S(e)=e$.  More
generally, for any $\gamma\in\GH$ denote
\begin{equation*}
L_\gamma := \{l\in H\mid al =\gamma(a) l,\ \forall a\in H\}\,,
\end{equation*}
then $L\equiv L_\ep$ is the space of left integrals.

From now on we assume $\dim H<\infty$.
To prove $\dim L=1$ we will show below that $H$ is a
Frobenius algebra, implying $\dim L_\gamma =1$ for all $\gamma\in\GH$.
Indeed, let $\om:H\to \kk$ be non degenerate and denote $\om_R:H\to\hat H,\ 
\bra\om_R(a)\mid b\ket:=\om(ba)$, then
\begin{equation*}
%\label{3.1}
  L_\gamma=\om_R^{-1}(\kk\gamma)\,.
\end{equation*}
Moreover, in this case we have for all $l\in L_\gamma$ and all $a\in H$
\begin{equation*}
la=\tilde\gamma(a)l,
\end{equation*}
where $\tilde\gamma=\gamma\circ\theta_\om^{-1}$, and where $\theta_\om =
\om^{-1}_R\circ \om_L \, \in\Aut H$ denotes the modular (or Nakayama)
automorphism of $\om$ (i.e. solving $\om(ab)=\om(b\theta_\om(a)),\ \forall
a,b\in H$).  Note that since $\theta_\om$ is unique up to inner automorphisms,
$\tilde\gamma$ only depends on $\gamma$.  As in ordinary f.d.
Hopf algebras, we call
\begin{equation}\label{3.2}
\mu:=\tilde\ep\equiv\ep\circ\theta_\om^{-1}\in\GH
\end{equation}
the {\em modulus} of $H$ ($\mu$ is also called the {\em distinguished
  grouplike element of $\hat H$} \cite{Rad93}).  Thus $H$ is unimodular, i.e.
$\mu = \ep$, iff one (and therefore all) nonzero integrals are two-sided. In
particular if $H$ is symmetric (i.e. admits a non degenerate trace), then it
is unimodular.

We now generalize the methods of \cite{LS} to prove that indeed all f.d.
quasi-Hopf algebras are Frobenius, such that the above arguments apply.
First, we consider $(H,\Delta)$ as a left quasi-Hopf $H$-bimodule and choose
the dual (right) quasi-Hopf $H$-bimodule structure  $(H^*,\rho)$ as in
Proposition~\ref{prop2.12}.  Thus, as a linear space $H^*=\dH$ with $H$-bimodule structure given for
$a,b\in H$ and $\psi\in H^*$ by (see \eqref{1.30} for the notation)
\begin{equation}\label{3.3'}
a\cdot\psi\cdot b= S(b)\arr\psi\arl S^{-1}(a)\,.
\end{equation}
Following Proposition~\ref{prop2.12} the $H$-coaction $\rho:H^*\to H^*\tp
H$ is given by 
\begin{equation}\label{3.3}
\rho(\psi):=\sum_i b^i*\psi\tp b_i\,,
\end{equation}
where $\{b_i\}\subset H$ is a $\kk$-basis with dual basis $\{b^i\}\subset\hat
H$, and where according to \eqref{2.16} the (non-associative) "multiplication" 
$*:H^*\tp H^*\to H^*$ is given by 
\begin{equation}\label{3.3''}
\bra \varphi * \psi\mid a\ket :=\bra \varphi\tp\psi\mid \bcop(a)\ket 
= \bra \varphi \tp \psi \mid V \cop(a) U \ket,\quad a\in H,\,\varphi,\psi\in H^*.
\end{equation}
Now in ordinary
Hopf algebras the coinvariants of $H^*$ would precisely be the left
integrals in $\hat H$ \cite{LS}.  Thus we propose
\begin{definition}
\label{def3.5} 
The coinvariants $\lambda\in \hat{H}^{coH}$ are called {\em left cointegrals
  on} $H$ and we denote the space of left cointegrals as $\L:=\hat{H}^{coH}\,.$
\end{definition}
We will give some
equivalent characterizations of left cointegrals in Section~5 and 
also in Section~7.
Theorem~\ref{thm2.7} and Proposition~\ref{prop2.12} now immediately imply
\begin{theorem}
\label{thm3.6} 
Let $H$ be a f.d. quasi-Hopf algebra. Then $\dim \L=1$ and all nonzero left
cointegrals on $H$ are nondegenerate. In particular $H$ is a Frobenius
algebra and therefore the space of left (right) integrals in $H$ is
one dimensional.
\end{theorem}
\begin{proof}
  The first statement follows from $\dim H=\dim H^*$ and the fact that by
  Theorem~\ref{thm2.7}
$\, \L\tp H\ni \lambda\tp a\mapsto \lambda\cdot a\equiv
  (S(a)\arr\lambda)\in H^* $
 provides an isomorphism of quasi-Hopf $H$-bimodules. For the second
 statement see the remarks above.
\end{proof}
Note that according to \eqref{2.3},\eqref{2.15} and
\eqref{3.3} the projection $E:H^*\to \L$ is given for $\varphi\in H^*$
and $a\in H$ by
\begin{equation}
\label{3.0}
 \bra E(\varphi)\mid a\ket =\sum_i \bra b^i\tp\varphi\mid
\bcop\big(S^{-1}(q^1_R)a  S^2(q_R^2 b_i) S(\beta )\big)\ket.
\end{equation}
\begin{lemma}
\label{lem3.7}
The transpose  $E^T : H \to H$  is given by
\begin{equation*}
E^T(a)= \sum_i ( b^i\tp\id)
\big(\bcop\big(S^{-1}(q^1_R)a  S^2(q_R^2 b_i) S(\beta )\big)\big)
\end{equation*}
and provides a projection onto the space of right integrals
$R\subset H$.  Moreover, the dual pairing 
$\L\tp R\ni\lambda\tp r\mapsto\bra\lambda\mid r\ket\in \kk$ is nondegenerate.
\end{lemma}
\begin{proof}
%  Using $b^i\tp xb_i=b^i\arl x\tp b_i$ we compute
%\begin{align*}
%  E^T(a) & = \bra b^i\tp\id\mid(S^{-1}(\alpha
%  Z^i)Y^i\tp\one)V\Delta(S^{-1}(X^i))\Delta(aS^2(b_i)S(\beta))U\ket
%  \\
%  & =\bra b^i\tp\id\mid h\Delta(aS^2(b_i)S(\beta))U\ket
%\end{align*}
%where we have used \eqref{2.18} and \eqref{?}. 
  By Lemma~\ref{lem2.6} and the definition \eqref{2.15} $E^T(H)=\{r\in H\mid
  rS(a)=r\ep(a),\ \forall a\in H\}$, which is the space of right integrals
  $R\subset H$. Thus $E$ and $E^T$ induce the splittings $H^*=\L\oplus
  R^\perp$ and $H=R\oplus\L^\perp$, respectively.
%%%%%%%%%%%%%%%%%%%%%%%%%%%%%%%%%%%
% Das gilt für jede Projektion E mit \L=Im E und R=Im E^T !
%%%%%%%%%%%%%%%%%%%%%%%%%%%%%%%%%%
%(since $\bra E(\varphi)\mid
%  a\ket = \bra E(\varphi)\mid E^T (a)\ket$, i.e. $\bra E(\varphi)\mid a\ket =
%  0 \quad\forall a\in H \Leftrightarrow \bra E(\varphi)\mid r \ket = 0$).
%%%%%%%%%%%%%%%%%%%%%%%%%%%%%%%%%%%
%  Überflüssig !
%%%%%%%%%%%%%%%%%%%%%%%%%%%%%%%%%%%
\end{proof}
Whether the projection onto the space of left integrals given by \cite{PO}  is
functorially related to our formula remains unclear at
the moment.

\section{Fourier Transformations}
In this Section we first determine the modular automorphism of a nonzero left cointegral in
terms of the modulus $\mu$ of $H$, just as for ordinary Hopf algebras.  
In particular, $H$ will be unimodular iff the
modular automorphism of any left cointegral on $H$ is given by the square of
the antipode. 
We then develop a notion of Fourier transformation for quasi-Hopf algebras and show
 that these are given in terms of cointegrals by the same formula as for ordinary Hopf algebras. 
This will finally lead to a characterization of symmetric or semi-simple f.d. quasi-Hopf algebras
just like in the coassociative case.

First note that by \eqref{2.4} and
Proposition~\ref{prop2.3} $\L$ carries a unital representation of $H$ given by
\begin{equation*}
a\re\lambda:=E(\lambda\arl S^{-1}(a)),\quad a\in H,\,\lambda\in \L\,.
\end{equation*}
Since $\dim \L=1$ there must exist a unique $\ga\in\GH$ such that
\begin{equation}
\label{3.1}
a\re\lambda=\ga(a)\lambda,\quad\forall a\in H,\,\lambda\in \L\,.
\end{equation}
Proposition~\ref{prop2.3} (v) then implies for all $a\in H,\,\lambda\in \L$
\begin{equation}
\label{3.4}
\lambda\arl S^{-1}(a)=\ga(a\1)\big(S(a\2)\arr\lambda\big) =
S(a\arl\ga)\arr \lambda.
\end{equation}
\begin{lemma}
\label{lem3.8} 
Let $H$ be a f.d. quasi-Hopf algebra with modulus $\mu$. Then $\ga=\mu$ and
the modular automorphism of any nonzero $\lambda\in \L$ is given by
\begin{equation*}
\theta_\lambda(a)=S\big(S(a)\arl\mu\big),\quad a\in H\,.
\end{equation*}
In particular, $H$ is unimodular 
 iff  $\lambda (ab) = \lambda (b S^2 (a))$ for all $\lambda \in \L$ and all $a,b\in H$.
\end{lemma}
\begin{proof}
By the defining relation of $\theta_\lambda$,
\begin{equation*}
    \lambda \arl a = \theta_\lambda(a) \arr \lambda,\quad \forall a \in H,
\end{equation*}
we conclude from \eqref{3.4}
$\theta_\lambda(a)=S\big(S(a)\arl\ga\big).$
This implies $\ga^{-1}\circ S^{-1}\circ\theta_\lambda=\ep\circ S=\ep$
and therefore $\mu\equiv\ep\circ\theta_\lambda^{-1}=\ga^{-1}\circ S^{-1}=\ga$.
\end{proof}
\begin{corollary}\label{cor5.2}
Let $r\in H$ be a right integral. Then $ar=\mu^{-1}(a)r$ for all $a\in H$.
\end{corollary}
\begin{proof}
Clearly $ar$ is a right integral and for any
$0\neq\lambda\in\L$ we have
$\bra\lambda\mid ar\ket=\bra S(a)\re\lambda\mid r\ket =\mu^{-1}(a)\bra\lambda\mid r\ket$,
from which the statement follows by the nondegeneracy of the pairing $\L\tp\ R\to\kk$.
\end{proof}
We now generalize the notion of a Fourier transformation.
Since $\dim\L=1$, any nonzero $\lambda\in\L$ induces an identification
$\L\tp H\cong H_\mu$ as quasi-Hopf bimodules, see Corollary~\ref{cor2.3}.
Thus, by Theorem~\ref{thm2.7} $H_\mu\cong H^*$ as quasi-Hopf $H$-bimodules.
\begin{definition}
\label{fourier}
A {\em Fourier transformation} is a nonzero morphism of quasi-Hopf $H$-bimodules
$\F:H_\mu\to H^*$, i.e. by \eqref{3.3'}-\eqref{3.3''} a linear map
satisfying  for  
$a,b\in H$, $\psi\in H^*$ and $T_\mu\in H\tp H$ defined in
Corollary~\ref{cor2.3}
\begin{align}
  \label{F1}
  \F(ab) & =S(b)\arr\F(a)=\F(b)\arl S^{-1}(a\arl\mu^{-1}) \\
  \label{F2}
  \psi*\F(a) & = (\F\tp\psi)(T_\mu\cop(a))\, .
\end{align}
\end{definition}
We will see that given \eqref{F2} the two
conditions in \eqref{F1} are actually equivalent. 
To this end, for $\lambda\in H^*$ define $\Fl:H\to H^*$ and $\Fl':H\to H^*$ by
\begin{equation}\label{F3}
\Fl(a):=S(a)\arr\lambda,\qquad \Fl'(a):=\lambda\arl S^{-1}(a\arl\mu^{-1})\,.
\end{equation}
Then by \eqref{F1} any Fourier transformation $\F$ satisfies $\F=\Fl=\Fl'$, where $\lambda:=\F(\one)$.
Moreover, as for ordinary Hopf algebras, $\lambda\equiv\F(\one)$ is a left cointegral.
More precisely, generalizing results of \cite{Nill94} we have
\begin{proposition}
\label{prop3.8}
Let $H$ be a f.d. quasi-Hopf algebra with modulus $\mu$ and let $0\neq\lambda\in\dH$.
Then the following are equivalent
\begin{itemize}
\item[(i)]
$\lambda$ is a left cointegral
\item[(ii)]
$\Fl$ is a Fourier transformation
\item[(iii)]
$\rho(\Fl(a))=(\Fl\tp\id)(T_\mu\cop(a)),\quad\forall a\in H$
\item[(iv)]
$\rho(\Fl'(a))=(\Fl'\tp\id)(T_\mu\cop(a)),\quad\forall a\in H$
\end{itemize}
\end{proposition}
\begin{proof}
(i)$\Rightarrow$(ii):
By Theorem~\ref{thm2.7} $\L\tp H\ni\lambda\tp a\mapsto\Fl(a)\in H^*$ is an 
isomorphism of quasi-Hopf $H$-bimodules, implying $\Fl:H_\mu\to H^*$
to be a Fourier transformation.

(ii)$\Rightarrow$(iii+iv):
Holds by condition \eqref{F2} in Definition~\ref{fourier}.

(iii)$\Rightarrow$(i):
Pick $0\neq\chi\in\L$ and put $f:=\F_\chi^{-1}\circ\F_\lambda$.
Since $\F_\lambda$ and $\F_\chi$ are right $H$-module maps, we get
$f(a)=f(\one)a,\ \forall a\in H$.
Moreover, (iii) implies $(f\tp\id)(T_\mu\cop(a))=T_\mu\cop(f(a)),\ \forall a\in H$,
and in particular $(f(\one)\tp\one)T_\mu=T_\mu\cop(f(\one))$.
Applying $\ep\tp\id$ gives $f(\one)=\ep(f(\one))\one$ and therefore
$\lambda=\ep(f(\one))\chi$.

The implication (iv)$\Rightarrow$(i) follows similarly by considering 
$f':=(\F'_\chi)^{-1}\circ\F'_\lambda$ and noting $f'(a)=af'(\one)$.
\end{proof}
Next, we determine the Frobenius basis associated with a non-zero left
cointegral $\lambda$ on $H$.  By this we mean the unique solution $\sum_i
u_i\tp v_i\equiv Q_\lambda\in H\tp H$ of
\begin{equation*}
\sum_i \lambda(au_i)v_i = a = \sum_i u_i\lambda(v_ia),\quad\forall a\in H\,.
\end{equation*}
Note that for $\dim H<\infty$ and $\lambda\in \hat H$, $Q_\lambda\in H\tp H$
as above exists if and only if $\lambda$ is nondegenerate.  
(Choose a basis $\{u_i\}$ of $H$ and let $\{v_i\}$ be the unique
basis of $H$ such that $\lambda (v_i u_j) = 
\delta_{i,j}$.) 
In this case the
two conditions above are equivalent and we have
$\sum_i au_i\tp v_i =\sum_i u_i\tp v_ia\,$
for all $a\in H$.
Also note that 
\begin{equation}
  \label{3.5}
  \sum_i u_i \tp v_i = \sum_i \theta_\lambda (v_i) \tp u_i.
\end{equation}
If $H$ is a f.d. Hopf algebra and $\lambda\in\hat H$ is a nonzero left
integral in $\hat H$, then by results of \cite{LS}
$Q_\lambda=(S\tp\id)(\Delta(r))$, where $r\in H$ is the unique right integral
satisfying $\bra\lambda\mid r\ket=1$.  For f.d.  quasi-Hopf algebras, this
result generalizes as follows:
\begin{proposition}
\label{prop3.9}
Let $H$ be a f.d. quasi-Hopf algebra and let $\lambda\in \L$ be a nonzero left
cointegral on $H$. Then
\begin{equation*}
\F_\lambda^{-1}(\psi)=(\id\tp\psi)(\bcop(r))\quad \mbox{and}\quad
Q_\lambda=(S\tp\id)(\bcop(r)),
\end{equation*}
where $r\in R\subset H$ is the unique right integral satisfying
$\bra\lambda\mid r\ket=1$.
\end{proposition}
\begin{proof}
  Since $\F_\lambda(a)=\nu(\lambda\tp a)$,  
  according to Theorem~\ref{thm2.7} and Lemma~\ref{lem3.7}
  the inverse $\F_\lambda^{-1}$ is given by
  \begin{equation*}
  \F_\lambda^{-1}(\psi)=\bra E(\psi\0)\mid r\ket \psi\1=\bra\psi\0\mid
  r\ket\psi\1=(\id\tp\psi)(\bcop(r))\,.
  \end{equation*}
  Here we have used $E^T(r)=r$ and the definition \eqref{2.21} applied to
  $(M,\Lambda):=(H,\Delta)$.  Hence we conclude $(\id\tp
  S(a)\arr\lambda)(\bcop(r))=a,\ \forall a\in H$, and therefore $\sum_i u_i\tp
  v_i:=(S\tp\id)(\bcop(r))$ satifies $\sum_i u_i\lambda(v_i a)=a,\ \forall
  a\in H$.
\end{proof}
We conclude this section by characterizing symmetric or semi-simple
f.d. quasi-Hopf algebras quite analogously as in the Hopf case.
\begin{proposition}
\label{prop3.10} 
A f.d. quasi-Hopf algebra $H$ is
\begin{itemize}
\item[(i)] unimodular if and only if one (and hence all) nonzero left (right)
integrals are $S$-invariant;
\item[(ii)] symmetric if and only if it is unimodular and $S^2$ is an inner
automorphism.
\end{itemize}
\end{proposition}
\begin{proof}
  (i): If there is a nonzero $S$-invariant integral in $H$,
  then it is two-sided, hence all integrals are two-sided and $S$-invariant
  (by the one-dimensionality of the space of left/right integrals), implying
  $\mu=\ep$.  Conversely, if $H$ is unimodular then all integrals are
  two-sided and by Lemma~\ref{lem3.8} the modular automorphism of any nonzero
  $\lambda\in \L$ is given by $S^2$. Thus by \eqref{3.5} the Frobenius basis
  $\sum_i u_i\tp v_i\equiv Q_\lambda\in H\tp H$ satisfies $\sum_i u_i\tp v_i=
  \sum_j S^2(v_j)\tp u_j$ and therefore $\sum_i \ep(u_i)v_i=\sum_j
  u_j\ep(v_j)$.  But with $(\ep\tp\id)(\bcop(r))=r$ and $(\id\tp\ep)(\bcop(r))=S^{-1}(\beta)r\alpha = r$
  (since $r$ is a two-sided integral and $\ep(S^{-1}(\beta)\alpha) = \ep
  (\beta \alpha) = 1$) Proposition~\ref{prop3.9}
  implies $r=\ep(u_i)v_i=u_i\ep(v_i)=S(r)$.

 (ii): If $H$ is symmetric then it is unimodular and all modular
automorphisms are inner. Hence $S^2$ is inner.  Conversely, if $H$ is
unimodular and $S^2$ is inner pick $0\neq\lambda\in \L$ and $g\in H$
invertible such that $gxg^{-1}=S^2(x),\ \forall x\in H$. Then
$\tau:=g^{-1}\arr\lambda$ is a nondegenerate trace on $H$.
\end{proof}
In \cite{P3} F. Panaite has shown recently, that a quasi-Hopf algebra is
semi-simple Artinian if and only if it contains a normalized left (or right)
integral. For finite dimensional quasi-Hopf algebras we have in addition (the
equivalence (v) in Theorem~\ref{thm3.11} is essentially also due to \cite{P3})
\begin{theorem}
\label{thm3.11} 
For a f.d. quasi-Hopf algebra $H$ the following are equivalent:
\begin{itemize}
\item[(i)] $H$ is semi-simple.
\item[(ii)] $H$ has a normalized left (or right) integral.
\item[(iii)] The left cointegral $\lambda_e:= E(\ep)\equiv
\sum_i S^2(b_i)S(\beta)\alpha\arr b^i\in\L$ is nonzero\footnote{where as above
  $b_i$ denotes a basis of $H$ with dual basis $b^i \in \dH$}.
\item[(iv)] $H$ has a Haar integral $e$.
\item[(v)] $H$ is a separable $\kk$-algebra.
\end{itemize}
Moreover, in this case the Haar integral $e$ satisfies $\bra \lambda_e\mid
e\ket =1$.
\end{theorem}
\begin{proof}
  (i)$\Rightarrow$(ii):  see Theorem 2.3 in \cite{P3}.  

  (ii)$\Rightarrow$(iii): Let $r\in R$ be a normalized right
  integral, then $\bra E(\ep)\mid r\ket = \ep(r)=1$. Moreover, to verify the
  above formula for $\lambda_e$, we use \eqref{3.0} and $(\ep\tp
  \id)(\bcop(x)) = S^{-1}(\beta)x\alpha$ to obtain
  \begin{align*}
    \lambda_e (a) &= \sum_i \bra b^i \mid S^{-1}(q_R^1 \beta) aS^2(\qr^2 b_i)
    S(\beta) \alpha \ket \\
     &= \sum_i \bra b^i \mid \qr^2 S^{-1}(q_R^1 \beta) aS^2( b_i)
    S(\beta) \alpha \ket\\
    &= \sum_i \bra b^i \mid  a S^2( b_i) S(\beta) \alpha \ket\\
  \end{align*}
and therefore indeed $\lambda_e= \sum S^2(b_i)S(\beta)\alpha\arr
  b^i$. Here we have 
used  \eqref{1.8} and the identity  $\sum b^i \tp y b_i = \sum b^i \arl y \tp b_i$.

(iii)$\Rightarrow$(iv): If $\lambda_e\neq 0$ there exists a unique right
integral $e\in H$ such that $\lambda_e(e) =1$ yielding $\ep(e)=\ep(E^T(e))=1$.
Thus $e$ is normalized.  Moreover, using again the identity $\sum b^i \tp y
b_i = \sum b^i \arl y \tp b_i$ one easily verifies that $\lambda_e (a S^2(b))
= \lambda(ba)$, i.e.  the modular automorphism of $\lambda_e$ is given by
$S^2$. Hence, by Lemma~\ref{lem3.8}, $H$ is unimodular and $ e$ is the
Haar integral in $H$. 

(iv)$\Rightarrow$(v): Following \cite{P3}, if $e$ is a normalized left
integral then
\begin{equation*}
P:=(\id\tp S)(q_R\Delta(e)(\beta\tp \one))
\end{equation*}
provides a separating idempotent in $H\tp H_{op}$, i.e. $P_1^iP_2^i=\one$
and $(a\tp\one)P=P(\one\tp a)$ for all $a\in H$.

(v)$\Rightarrow$(i): This is a standard textbook exercise, see
  e.g. \cite{Pierce}. 
\end{proof}

Let us conclude with an explicit formula for $\lambda_e$ in the case of $H$
being a multi-matrix algebra, i.e.  $H\cong\oplus_I\Mat\!_\kk(n_I)\,$.  In this
case $H$ is symmetric, i.e. it admits a non degenerate trace.  Hence
Proposition~\ref{prop3.10} applies and we may pick $g\in H$ invertible such that
$gxg^{-1}=S^2(x),\ \forall x\in H$. Put $c_I:=tr_I(g^{-1}S(\beta)\alpha)$,
where $tr_I$ denotes the standard trace on $\Mat\!_\kk(n_I)$, and
denote $\LL(a)$, $\RR(b)$ the operators of left multiplication with $a$ and right
multiplication with $b$, respectively, on $H$. Then
\begin{align}
  \notag
  \lambda_e(a) &= \sum_i \bra b^i \mid a g b_i g^{-1} S(\beta) \alpha \ket \\
   &= Tr\big(\LL(ag)\circ\RR(g^{-1}S(\beta)\alpha)\big)\notag \\
\label{3.10}
&= \sum_I c_I tr_I(ag),
\end{align}
 where $Tr$ denotes the standard trace on
$\End\!_\kk H$ and where we have used that $Tr(\LL(x) \RR(y)) = \sum_I tr_I (x) tr_I (y)$. 
By the nondegeneracy of $\lambda_e$ we conclude $c_I\neq 0,\ \forall I$.
Proposition~\ref{prop3.9} then implies for the Haar integral $e\in H$
\begin{equation*}
(S\tp\id)(\bcop(e)) =\sum_{I,\mu,\nu}\frac{1}{c_I}e_I^{\mu\nu}g_I^{-1}\tp
e_I^{\nu\mu},
\end{equation*}
where the $e_I^{\mu\nu}$'s denote the matrix units in $\Mat\!_\kk(n_I)$
and $g_I$ the component of $g$ in the matrix block $\Mat\!_\kk(n_I)$.

%%% Local Variables: 
%%% mode: latex
%%% TeX-master: "integral"
%%% End: 

\section{The Comodulus and Radford's Formula}
Recall that if $H$ is a finite dimensional Hopf algebra we may identify the modulus of $\dH$
with an element $u\in H$, i.e.
\begin{equation*}
\lambda \psi = \psi(u) \,\lambda
\end{equation*}
for all left cointegrals $\lambda\in\L$ and all $\psi\in\dH$.
If $r\in H$ satisfies $\lambda(r)=1$ we get
\begin{equation*}
u=(\lambda\tp\id)\big(\cop(r)\big)\, .
\end{equation*}
Choosing $r$ to be a right integral, $r\in R$, this is the definition which appropriately
generalizes to quasi-Hopf algebras:
\begin{definition}
\label{def5.1}
Let $H$ be a f.d. quasi-Hopf algebra and let $\lambda \in \L$ and $r\in R$ 
satisfy $\bra\lambda\mid r\ket =1$.
Then we call
\begin{equation*}
  u:=(\lambda\tp\id)(\bcop(r))\equiv(\lambda\tp\id)(V\cop(r)U)\in H
\end{equation*}
the {\em comodulus} of $H$.
\end{definition}
The good use of this definition stems from the fact that it gives rise to a generalization
of Radfords Formula \cite{Rad76} expressing $S^4$ 
as a composition of the inner and coinner automorphisms, respectively,
induced by $u^{-1}$ and $\mu$.
\begin{proposition}\label{prop5.2} {\rm (Radford's Formula)}
Let $H$ be a f.d. quasi-Hopf algebra with modulus $\mu\in\dH$ and comodulus $u\in H$.
For $a\in H$ put $S_\mu(a):=S(a)\arl\mu$ and $v:=(\id\tp\lambda\circ S)(\bcop(r))$.
Then
\begin{align}
  u^{-1}&=S^2(v)=S_\mu^{-2}(v) 
  \\
  \label{6.2}
  u^{-1}\,a\,u&=S^2(S_\mu^2(a)),\quad\forall a\in H\,.
\end{align}
\end{proposition}
\begin{proof}
Using \eqref{2.19} and Corollary~\ref{cor5.2} we have for all $a,b\in H$
\begin{equation*}
 [\one\tp S^{-1}(a)]\,\bcop(r)\,[\one\tp S(b)] =
 [a\arl\mu\tp\one]\,\bcop(r)\,[b\tp\one] 
\end{equation*}
and therefore
\begin{align*}
  au &=\bra\lambda\tp\id\mid[S_\mu(a)\tp\one]\bcop(r)\ket
  \\
  &=\bra\lambda\tp\id\mid\bcop(r)[S(S^2_\mu(a))\tp\one]\ket
  \\
  &=uS^2(S^2_\mu(a)),
\end{align*}
where we have used that according to Lemma~\ref{lem3.8} the modular automorphism 
of $\lambda$ is given by $\theta_\lambda=S\circ S_\mu$.
We are left to show that $S_\mu^{-2}(v)$ is a left inverse of $u$, implying by the above calculation
$S^2(v)$ to be a right inverse and therefore $u^{-1}=S^2(v)=S_\mu^{-2}(v) $.
First note that the second line of the above calculation also gives
\begin{equation*}
  au=\bra\F_\lambda(S^2_\mu(a))\tp\id\mid\bcop(r)\ket
\end{equation*}
for all $a\in H$.
Hence
\begin{align*}
  S^{-2}_\mu(v)u &=\bra\F_\lambda(v)\tp\id\mid\bcop(r)\ket
  \\
  &=(\lambda\tp\id)(Q_\lambda) =\one
\end{align*}
where we have used that by Proposition~\ref{prop3.9}
$v=\F_\lambda^{-1}(\lambda\circ S)$ and that $Q_\lambda:=(S\tp\id)(\bcop(r))\in H\tp H$
provides the Frobenius basis of $\lambda$.
\end{proof}
%Now $ \mu \circ S = \mu^{-1}$ implies
%$S(a) \arl \mu = S(\mu^{-1} \arr a )$ and 
%$\mu^{-1} \arr S(a) = S(a \arl \mu)$. Thus we get  
%$
%  S^2_\mu (a) = S^2\big( (\mu^{-1}\arr a) \arl \mu \big)$ 
%and therefore \eqref{6.2} yields
%\begin{equation*}
%  S^4 (a) = u^{-1} \, \big( (\mu \arr a ) \arl \mu^{-1}\big) \,u\, .
%\end{equation*}
To get more similarity with Radford's original formulation in \cite{Rad76} we put
$f_\mu:=(\mu\tp\id)(f)$ and use \eqref{1.12} to get 
$S^2_\mu(a)=f_\mu^{-1}\,S(\mu^{-1}\arr(S(a)\arl\mu))\,f_\mu$. Writing 
$b=\mu^{-1}\arr(S(a)\arl\mu)$ we conclude
\begin{corollary}
Under the conditions of Proposition~\ref{prop5.2} we have for all $b\in H$
$$
S^4(b) = S^3(f_\mu^{-1})S(u)[(\mu\arr b)\arl\mu^{-1}]S(u^{-1})S^3(f_\mu)
$$
\end{corollary}

\section{Cocentral Bilinear Forms} 
Recall \cite{Abe} that a Hopf algebra
$H$ is cosemisimple iff it admits a bilinear form $\Sigma: H\tp H\to \kk$
satisfying for all $a,b\in H$
\begin{align}
\label{4.1}
(\id\tp\Sigma)\big(\Delta(a)\tp b\big) &=
(\Sigma\tp\id)\big(a\tp\Delta(b)\big), \\ 
\label{4.2} 
\Sigma\circ\Delta &= \ep.
\end{align}
In particular, if $\lambda\in\hat H$ is a normalized left cointegral on $H$
(i.e. $\lambda(\one)=1$ and $\lambda\arr a=\lambda(a)\one,\ \forall a\in H$),
then such a $\Sigma$ is obtained by
\begin{equation*}
\Sigma(a\tp b):=\lambda(aS(b))\,.
\end{equation*}
In this case $\Sigma$ is  {\em right invariant}, i.e. for all $a,b,c\in H$
\begin{equation*}
\Sigma\big((a\tp b)\cop(c)\big) = \Sigma(a\tp b)\ep(c)\,.
\end{equation*}
Conversely, any right invariant bilinear form $\Sigma$ is obtained this
way, where $\lambda(a):=\Sigma(a\tp \one)$, and where $\Sigma$ satisfies
\eqref{4.1} iff $\lambda\in \hat H$ is a left cointegral on $H$. Moreover, in
this case the normalization condition \eqref{4.2} is equivalent to
$\Sigma(\one\tp\one)\equiv\lambda(\one)=1$.

In this Section we generalize these relations to f.d. quasi-Hopf
algebras.  Our motivation is two-fold.  On the one hand, this will provide a
characterization of left cointegrals which is less implicit than
Definition~\ref{def3.5} and more reminiscent to ordinary Hopf algebra theory.
On the other hand, although we have no sensible notion of cosemisimplicity for
quasi-Hopf algebras, we precisely need the analogues of \eqref{4.1} and
\eqref{4.2} to prove in Theorem~\ref{thm4.8} semisimplicity of
the diagonal crossed 
products $A\reli \hat H$ constructed in \cite{HN1}, for any semisimple
two-sided $H$-comodule algebra $A$.
In particular this will imply the quantum double $D(H)$ of a f.d.
quasi-Hopf algebra $H$ to be semisimple, iff $H$ is semisimple and admits a
  left cointegral $\lambda_0\in \L$ satisfying the normalization condition
  $\lambda_0(\beta S(\alpha))=1$.    

We start with first preparing a more general formalism

\begin{definition}
\label{6.1}
Let $(K,\Lambda_K)$ and $(M,\rho_M)$ be a left and a right quasi-Hopf $H$-bimodule, respectively.
We call a bilinear form $\Sigma: K\tp M\to\kk$ {\em $H$-biinvariant}, if for all $a \in H,\ k\in K,\ m\in M$
\begin{equation*}
\Sigma\big((k\tp m)\cdot\cop(a)\big) = \ep(a)\Sigma(k\tp m) = 
  \Sigma\big(\cop(a)\cdot(k\tp m)\big)\,.
\end{equation*}
We call $\Sigma$ {\em $H$-cocentral}, if 
\begin{equation*}
(\id\tp\Sigma)\big(\phi\cdot[\Lambda_K(k)\tp m]\cdot\phi^{-1}\big) =
(\Sigma\tp\id)\big(\phi^{-1}\cdot[k\tp\rho_M(m)]\cdot\phi\big)
\end{equation*}
\end{definition}
If the $H$-(co)actions are understood we also just use the words {\em biinvariant} and {\em cocentral}.
 To motivate our terminology note that if H is a Hopf algebra, then $(K\tp M)^\wedge$ naturally becomes an 
 $\dH$-bimodule by putting 
$$
(\psi\cdot\Sigma\cdot\varphi)(k\tp m):= (\varphi\tp\Sigma\tp\psi)(\Lambda_K(k)\tp\rho_M(m)),\
\varphi,\psi\in\dH,\ \Sigma\in (K\tp M)^\wedge\,.
$$
In this case $\Sigma$ is $H$-cocentral iff $\varphi\cdot\Sigma=\Sigma\cdot\varphi,\ \forall\varphi\in\dH$.

Let now $K^*=\Hom\!_\kk(K,\kk)$ be the $H$-bimodule with $H$-actions defined in \eqref{2.15}.
Recall from Proposition~\ref{prop2.12} that if $H$ or $K$ is finite dimensional, then there exists a right 
$H$-coaction $\rho_{K^*}:K^*\to K^*\tp H$ making $(K^*,\rho_{K^*})$ a right quasi-Hopf $H$-bimodule.
We are aiming to characterize cocentral biinvariant elements $\Sigma\in(K\tp M)^\wedge$
in terms of quasi-Hopf $H$-bimodule morphisms $f:M\to K^*$.
To this end we need
 \begin{lemma}
\label{lem4.5} 
Let $(K,\Lambda)$ be a left quasi-Hopf $H$-bimodule and (assuming $K$ or $H$ finite dimensional)
let $(K^*,\rho)$ denote its dual right quasi-Hopf $H$-bimodule, see Proposition~\ref{prop2.12}.
Then for all $k\in K$ and $\psi\in K^*$
\begin{equation*}
(\id\tp\psi) \big(q_R\cdot\Lambda(k)\cdot p_R\big) = ( k\tp\id )
\big( q_L\cdot\rho(\psi)\cdot p_L \big).
\end{equation*}
\end{lemma}
\begin{proof}
We first show the identities
\begin{align}
  \label{h1}
  q_R & = [q_L^2\tp\one]\,V\,\Delta\big(S^{-1}(q_L^1)\big),
  \\
   \label{h2}
  p_R & = \Delta\big(S(p_L^1)\big)\,U\,[p_L^2\tp\one].
\end{align}
Noting  that in $H_{op}$ the roles of $U$ and $V$ as well as of $q_R$
and $p_R$ and of $q_L$ and $p_L$ interchange reduces \eqref{h2} to \eqref{h1}.
For the  left hand side of \eqref{h1} we obtain, using
\eqref{2.18},\eqref{1.14} 
\begin{align}
  \notag
  [q_L^2\tp\one]\,V\,&\Delta\big(S^{-1}(q_L^1)\big) 
     \overset{\eqref{1.12}}{=}
        (S^{-1}\tp S^{-1})\big(f^{21} \, \cop^{op}(q^1_L) \, p^{21}_R \,
     [S(q_L^2) \tp\e]\big) \\\notag
    & \overset{\eqref{1.13}}{=} 
      (S^{-1}\tp S^{-1})\Big((S\tp S)(\cop(\Xbar^i))\,\gamma^{21} \,
    \cop^{op}(\Ybar^i) \, p^{21}_R \,  [S(\Zbar^i) \tp\e]\Big) \\\notag
    & \overset{\eqref{1.9'}}{=} 
      (S^{-1}\tp S^{-1})\Big((S\tp S)\big(X^j \tp Y^j \,
    \cop(\Xbar^k\Xbar^i)\big)\,[\alpha \tp \alpha]\,\\\notag
     &\qquad\, [\Zbar^k \tp Z^j\Ybar^k]\, \cop^{op}(\Ybar^i) \,[\Ybar^m \tp
    \Xbar^m]\, [\beta S(\Zbar^i\Zbar^m)\tp \e]\Big) \\\label{h3}
    &\quad = (S^{-1}\tp S^{-1}) \circ \sigma (\Psi),
\end{align}
where $\sigma : H^{\tp^5} \to H^{\tp^2}$ denotes the map
\begin{equation*}
  a\tp b \tp c \tp d\tp e \mapsto S(a)\alpha d \beta S(e) \tp S(b) \alpha c,
\end{equation*}
and $\Psi \in H^{\tp^5} $ is given by
\begin{equation*}
  \Psi := [\phi\tp\e\tp\e]\, [(\cop\tp\id\tp\id)(\phi^{-1})\tp\e] \,
         (\cop\tp\cop\tp\id)(\phi^{-1})\, [\e\tp\e\tp \phi^{-1}].
\end{equation*}
Now we use the pentagon equation for the last three factors of $\Psi$ to obtain
\begin{align*}
  \Psi &= [\phi\tp\e\tp\e]\,\big((\cop\tp\id)
  \cop \tp\id\tp\id\big)(\phi^{-1})  \, (\cop\tp\id\tp\cop)(\phi^{-1})\\
  &= \big((\id\tp\cop)\cop \tp\id\tp\id\big)(\phi^{-1})
  [\phi\tp\e\tp\e]\, (\cop\tp\id\tp\cop)(\phi^{-1}) .
\end{align*}
Using the antipode properties \eqref{1.7}, under the evaluation of
$\sigma$ the third factor may be dropped and  the first factor
may be replaced by $(\phi^{-1})^{145}$.
Hence we get
\begin{align*}
  \sigma (\Psi) &= \sigma\big( [\Xbar^i \tp\e\tp\e\tp \Ybar^i\tp \Zbar^i]\,
                     [\phi\tp\e\tp\e] \big) \\
     &= S(X^k)S(\Xbar^i)\alpha \Ybar^i \beta S(\Zbar^i) \tp S(Y^k)\alpha Z^k\\
    &= S(X^k) \tp S(Y^k)\alpha Z^k,
\end{align*}
where we have used \eqref{1.8}. Thus we finally arrive at
\begin{equation*}
  (S^{-1}\tp S^{-1})\big(\sigma (\Psi)\big) = X^k\tp S^{-1}(\alpha Z^k) Y^k =
  q^R .
\end{equation*}
By \eqref{h3} we have proved \eqref{h1}.

Using \eqref{h1},\eqref{h2} and the definition of $\rho$ given in 
\eqref{2.15}, \eqref{2.16} and  \eqref{2.21}, one obtains 
\begin{align*}
( k\tp\id ) \big( q_L\cdot\rho(\psi)\cdot p_L\big)  & = 
\big( S^{-1}(q^1_L)\cdot k\cdot S(p^1_L) \tp\id \big)
\Big([\e \tp q^2_L] \cdot\rho(\psi)\cdot [\e\tp p^2_L]\Big) 
\\
& = (\id\tp\psi)\Big( [q^2_L \tp \e]\cdot \bar\Lambda\big( S^{-1}(q^1_L)\cdot k\cdot
S(p^1_L)\big)\cdot [p^2_L\tp \e]\Big) 
\\
& = (\id \tp\psi)\big( q_R\cdot\Lambda(k)\cdot p_R \big).
\end{align*}
This finishes the proof of Lemma~\ref{lem4.5}.
\end{proof}
\begin{theorem}
\label{thm6.2}
Let $(K,\Lambda_K)$ and $(M,\rho_M)$ be a left and a right quasi-Hopf $H$-bimodule, respectively.
\begin{itemize}
 \item[(i)] Then the assignment $f\mapsto\Sigma_f$, where
\begin{equation*}
\Sigma_f(k\tp m) := \bra f(m)\mid S^{-1}(\alpha)\cdot k\cdot\beta\ket\,,
\end{equation*}
provides a one-to-one correspondence between $H$-bimodule maps $f:M\to K^*$
and biinvariant elements $\Sigma_f\in(K\tp M)^\wedge$. The inverse assignement
is given by $\Sigma\mapsto f_\Sigma$, where
\begin{equation*}
\bra f_\Sigma(m)\mid k\ket := \Sigma(p\cdot(k\tp m)\cdot q)
\end{equation*}
and where $p=p_L$ or $p=p_R$ and $q=q_L$ or $q=q_R$, 
$f_\Sigma$ being independent of these choices.
\item[(ii)] Assume in addition $H$ or $K$ finite dimensional. Then
$\Sigma_f$ is cocentral 
if and only if $f$ is a morphism of quasi-Hopf $H$-bimodules, 
i.e. $(f\tp\id)\circ \rho_M = \rho_{K^*}\circ f$.
\end{itemize}
\end{theorem}
\begin{proof}
For $a,b\in H$ we have by \eqref{1.7}
  \begin{align*}
\Sigma_f\big(\cop(a)\cdot(k\tp m)\cdot\cop(b)\big) & =
\bra a\2\cdot f(m)\cdot b\2\mid S^{-1}(\alpha)a\1\cdot k\cdot b\1\beta\ket
\\
&=\bra f(m)\mid S^{-1}(S(a\1)\alpha a\2)\cdot k\cdot b\1\beta S(b\2)\ket
\\
&=\ep(a)\ep(b)\Sigma_f(k\tp m)
\end{align*}
and from \eqref{1.8} we conclude  for $p=p_{L/R}$ and $q=q_{L/R}$
  \begin{align*}
\Sigma_f(p\cdot(k\tp m)\cdot q) &=
\bra p^2\cdot f(m)\cdot q^2\mid S^{-1}(\alpha)p^1\cdot k\cdot q^1\beta\ket
\\
&=\bra f(m)\mid S^{-1}(S(p^1)\alpha p^2)\cdot k\cdot q^1\beta S(q^2)\ket
\\
&=\bra f(m)\mid k\ket\,.
\end{align*}
Thus, $\Sigma_f$ is biinvariant and $f_{\Sigma_f}=f$ independently of the choice of
$q=q_{L/R}$ and $p=p_{L/R}$.
Conversely, let $\Sigma\in(M\tp K)^\wedge$ be biinvariant, then for all choices
$q=q_{L/R}$ and $p=p_{L/R}$
  \begin{align*}
\Sigma\big(p\cdot(k\tp a\cdot m)\big) &=\Sigma\big(p\cdot(S^{-1}(a)\cdot k\tp m)\big)
\\
\Sigma\big((k\tp m\cdot b)\cdot q\big) &=\Sigma\big((k\cdot S(b)\tp m)\cdot q\big)
\end{align*}
by the identities \eqref{1.22} - \eqref{1.23'}.
Hence, $f_\Sigma$ is an $H$-bimodule map for all choices of
$q=q_{L/R}$ and $p=p_{L/R}$:
  \begin{align*}
\bra f_\Sigma(a\cdot m\cdot b)\mid k\ket &=
\Sigma\big(p\cdot(k\tp a\cdot m\cdot b)\cdot q\big)
\\
&=\Sigma\big(p\cdot(S^{-1}(a)\cdot k\cdot S(b)\tp m)\cdot q\big)
\\
&=\bra a\cdot f_\Sigma(m)\cdot b\mid k\ket\,.
\end{align*}
Finally, if $\Sigma$ is biinvariant we may use the identities \eqref{1.24} and \eqref{1.25}
to conclude
  \begin{align*}
\Sigma\big((k\cdot\beta\tp m)\cdot q\big) =\Sigma\big((k\tp m\cdot S^{-1}(\beta))\cdot q\big)
&=\Sigma(k\tp m)
\\
\Sigma\big(p\cdot(S^{-1}(\alpha)\cdot k\tp m)\big) =\Sigma(p\cdot(k\tp\alpha\cdot m)\big)
&=\Sigma(k\tp m)
\end{align*}
for all choices of $q=q_{L/R}$ and $p=p_{L/R}$.
Hence, $\Sigma_{f_\Sigma}=\Sigma$ for all these choices and therefore $f_\Sigma$ is independent
of these choices. This proves the first statement of Theorem~\ref{thm6.2}.
To prove the second statement we compute
\begin{align*}
  (\id\tp\Sigma_f)&\big(\phi\cdot[\Lambda(k)\tp m]\cdot\phi^{-1}\big)\\
  & = \bra \id\tp Z^i\cdot f(m)\cdot\Zbar^j\mid 
  (X^i\tp S^{-1}(\alpha)Y^i)\cdot\Lambda(k)\cdot(\Xbar^j\tp \Ybar^j\beta)\ket
  \\
  & = \bra\id\tp f(m)\mid q_R\cdot\Lambda(k)\cdot p_R\ket
  \\
  & = \bra k\tp\id\mid q_L\cdot\rho_{K^*}(f(m))\cdot p_L\ket
\end{align*}
where in the last line we have used Lemma~\ref{lem4.5}.
On the other hand
\begin{align*}
  (\Sigma_f\tp\id) & \big(\phi^{-1}\cdot[k\tp\rho_M(m)]\cdot\phi\big) \\
  & = \bra f(\Ybar^i\cdot m\0\cdot Y^j)\mid S^{-1}(\alpha)\Xbar^i\cdot
  k\cdot X^j\beta\ket \,\Zbar^i m\1 Z^j
  \\
  & = \bra f(q_L^1\cdot m\0\cdot p_L^1)\mid k\ket\, q_L^2 m\1p_L^2\,.
\end{align*}
Thus, $\Sigma_f$ is cocentral if and only if for all $m\in M$
\begin{equation}
  \label{cocentral}
  q_L\cdot\rho_{K^*}(f(m))\cdot p_L = 
     (f\tp\id)\big(q_L\cdot\rho_M(m)\cdot p_L\big)\,.
\end{equation}
Clearly, \eqref{cocentral} holds, if $f$ is an $H$-bimodule and an $H$-comodule map.
Conversely, let $f$ be an $H$-bimodule map satisfying \eqref{cocentral}.
Then by \eqref{1.25}
\begin{align*}
 \rho_{K^*}(f(m))&=
 [(S(p_L^1)\tp\one)q_L\Delta(p_L^2)]\cdot\rho_{K^*}(f(m))\cdot
  [\Delta(q_L^2)p_L(S^{-1}(q_L^1)\tp\one)]
  \\
  &= [(S(p_L^1)\tp\one)q_L]\cdot\rho_{K^*}\big(f(p_L^2\cdot m\cdot q_L^2)\big)\cdot
  [p_L(S^{-1}(q_L^1)\tp\one)]
  \\
  &= [S(p_L^1)\tp\one]\cdot(f\tp\id)\big( q_L\cdot\rho_M(p_L^2\cdot m
  \cdot q_L^2)\cdot p_L \big) \cdot [ S^{-1}(q_L^1)\tp\one]
  \\
  &=
  (f\tp\id)\Big([(S(p_L^1)\tp\one)q_L\Delta(p_L^2)]\cdot
  \rho_M(m)\cdot[\Delta(q_L^2)p_L(S^{-1}(q_L^1)\tp\one)]\Big)  
  \\
  &= (f\tp\id)\big(\rho_M(m)\big)
\end{align*}
Hence $f$ is an $H$-comodule map. This proves the second statement of Theorem~\ref{thm6.2}.
\end{proof}
We now apply Theorem~\ref{thm6.2} to our theory of integrals by putting $(K,\Lambda_K)=(H,\Delta)$
and $(M,\rho_M)=(H_\mu,\rho_\mu)$, see Corollary~\ref{cor2.3}, 
where $\mu\in\Gamma(H)$ is the modulus of $H$.
In this case, by Definition~\ref{fourier} and Proposition~\ref{3.3} the quasi-Hopf bimodule
morphisms $H_\mu\to H^*$ are precisely our Fourier transformations $\F_\lambda,\ \lambda\in\L$,
given in \eqref{F3}.
Thus, we conclude
\begin{corollary}\label{cor6.4}
Let $H$ be a finite dimensional quasi-Hopf algebra with modulus $\mu$.
Then the assignment $\lambda\mapsto\Sigma_\lambda\equiv\Sigma_{\F_\lambda}$, i.e.
$$
\Sigma_\lambda(a\tp b):=\lambda(S^{-1}(\alpha)a\beta S(b))\,,
$$
provides a one-to-one correspondence between left cointegrals $\lambda\in\L$ and
cocentral biinvariant forms $\Sigma:H\tp H_\mu\to\kk$.
The inverse assignment is given by $\Sigma\mapsto\lambda_\Sigma:=f_\Sigma(\one)$, i.e.
(putting $p_\mu:=p^1\tp p^2\arl\mu$)
$$
\lambda_\Sigma(a):=\Sigma\left(p_\mu(a\tp\one)q\right)
$$
which is independent of the choices
$q=q_{L/R}$ and $p=p_{L/R}$.
\end{corollary}
Corollary~\ref{cor6.4} shows that on a finite dimensional quasi-Hopf algebra the space of 
biinvariant cocentral forms $\Sigma:H\tp H_\mu\to\kk$ is one dimensional.
Moreover, we also have
\begin{corollary}\label{cor6.5}
Let $H$ be a finite dimensional quasi-Hopf algebra with modulus $\mu$
and let $\gamma\in\Gamma(H)$.
If there exists a nonzero cocentral biinvariant form $\Sigma:H\tp H_\gamma\to\kk$,
then $\gamma=\mu$.
\end{corollary}
\begin{proof}
Pick $0\neq\lambda\in\L$ and put $f':=\F_\lambda^{-1}\circ f_\Sigma$. Then $f':H_\gamma\to H_\mu$
is a nonzero morphism of quasi-Hopf $H$-bimodules, whence
$f'(a)=ba=(a\arl\gamma^{-1}\mu)b$ for all $a\in H$, where $b:=f'(\one)\neq 0$.
Moreover, $T_\mu\cop(b)\equiv\rho_\mu(f'(\one))=(f'\tp\id)(T_\gamma)=(b\tp\one)T_\gamma$.
Applying $\ep\tp\id$ gives $b=\ep(b)\one$ and therefore $\gamma=\mu$.
\end{proof}
Putting $\gamma=\ep$ Corollary~\ref{cor6.5} in particular implies that nonzero cocentral
biinvariant forms $\Sigma:H\tp H\to\kk$ exist if and only if $H$ is unimodular.

Let us conclude with considering the normalization condition \eqref{4.2}.
We call a cocentral biinvariant $\Sigma:H\tp H_\mu\to\kk$ {\em normalized},
if $\Sigma(\one\tp\one)=1$, which by the right $H$-invariance is equivalent to \eqref{4.2}.
If such a normalized $\Sigma$ exists, it is unique and we denote it by $\Sigma_0$.
Correspondingly, we also call the associated left cointegral $\lambda_0$ normalized,
which means $\lambda_0(S^{-1}(\alpha)\beta)=1$.

It is instructive to look at the condition for the existence of $\lambda_0$
(and therefore $\Sigma_0$) in the case of $H$ being a multi-matrix algebra
over an algebraically closed field $\kk$.  In this case $\mu=\ep$ and we have to check, wether
the nonzero left cointegral $\lambda_e$ of Theorem~\ref{thm3.11} (iii) satisfies
$\lambda_e(S^{-1}(\alpha)\beta)\equiv\lambda_e(\beta S(\alpha))\neq 0$,
 which would give $\lambda_0=\lambda_e(\beta S(\alpha))^{-1}\,\lambda_e$.
If $\kk$ is algebraically closed we may
choose $g$ in \eqref{3.10} obeying 
\footnote{$S(g)g$ being a central
  $S$-invariant invertible element in $H$ it has a central $S$-invariant
  invertible square root $z\in Z(H)$. Thus, $g':=z^{-1}g$ satisfies
  $S(g')={g'}^{-1}$.}  
$S(g)=g^{-1}$ to conclude
\begin{equation}
\lambda_e(\beta S(\alpha))=\sum_I d_I d_{\bar I} =:d_H,
\end{equation}
where
\begin{equation*}
  d_I=tr_I(\beta S(\alpha) g),\quad  
  \quad d_{\bar I}=tr_I(S(\beta S(\alpha) g)).
\end{equation*}
Note that in the ribbon quasi-Hopf scenario of \cite{AC} the $d_I$'s are
precisely the quantum dimensions of
the irreducible $H$-modules labelled by $I$.
\footnote{In our description the $d_I$'s are not yet uniquely fixed, however
the products $d_Id_{\bar I}$ are so.}.  
Thus one might call $d_H$ the {\em quantum dimension of $H$},
considered as an $H$-module with respect to the 
adjoint action.  In summary we conclude
\begin{corollary}
\label{cor4.7} A f.d. semisimple quasi-Hopf algebra $H$ over an algebraically
closed field $\kk$ admits a normalized left cointegral $\lambda_0$
{\rm(}and therefore
a  normalized cocentral biinvariant form
$\Sigma_0:H\tp H\to \kk${\rm)}, if and only if the quantum dimension of $H$
with respect to the
adjoint action on itself is nonzero.
\end{corollary}

\section{Semisimplicity of Diagonal Crossed Products}
As an application we will now prove a Maschke type Theorem for
diagonal crossed products $A\reli \dH$ as defined in \cite{HN1,HN2}. 
We recall that by an twosided $H$-comodule algebra $(A,\delta,\Psi)$
associated with a quasi-Hopf-algebra
  $H$  we mean an algebra $A$ together with an algebra map 
  $\delta :\,A \to H\tp A\tp H$ and an invertible
  element $\Psi \in  H\tp H\tp A\tp H\tp H$ satisfying (denoting $\id
\equiv \id_H, \, \e \equiv \e_H$)
  \begin{subequations}
  \begin{align}
    \label{4.5a} 
    &(\id \tp\delta\tp\id )(\delta(a)) \,\Psi  = \Psi
    \,(\cop\tp\id_A \tp \cop)( \delta(a)), \quad \forall a \in A, \\ 
      \notag  (\e \tp\Psi\tp \e ) \,
   & (\id \tp\cop\tp\id_A\tp\cop\tp\id)(\Psi) \,
    (\phi\tp\e_A\tp\phi^{-1})
    \\ &=
    (\id\tp\id\tp\delta\tp\id\tp\id )(\Psi)\, 
    (\cop\tp\id \tp\id_A\tp\id\tp\cop)(\Psi), \label{4.5b} \\ 
   & (\ep\tp\id_A\tp\ep)\circ \delta = \id_A, \label{4.5c} \\ 
    (\id\tp\ep  \tp\id_A & \tp\ep\tp\id)(\Psi) =
    (\ep\tp\id\tp\id_A\tp\id\tp\ep)(\Psi) =
    \e \tp\e_A\tp\e. \label{4.5d} 
    \end{align} 
 \end{subequations}
We call $(\delta,\Psi)$ a two-sided $H$-coaction  on $A$ and denote
\begin{align}
 \notag
\delta (a)& = a_{(-1)} \tp a\0 \tp a\1\\
\notag
 \varphi\re a\li\psi &: =
(\psi\tp\id\!_A\tp\varphi)(\delta(a)), \quad a\in A, \, \varphi,\psi \in\dH, \\
  \label{4.6}
  \Om &:= (h^{-1})^{21}\, (\Si\tp\Si\tp\id_A\tp\id\tp\id)(\Psi),
\end{align}
with $h$ defined in \eqref{1.14}.
\begin{definition}\cite{HN1}
\label{def4.8}
Given a two-sided coaction $(\delta,\Psi)$ of a quasi-Hopf algebra $H$ on
a unital algebra $A$ we define the {\em right diagonal crossed product}
$A \reli\dH$ to be the vector space $A\tp\dH$ with multiplication rule
\begin{align*}
 & (a\tp \varphi)( b \tp\psi) := 
ab\0\Om^3\tp\big[(\Om^2S^{-1}(b_{(-1)})\arr\varphi\arl b\1\Om^4)(\Om^1\arr\psi\arl\Om^5) \big]
\\ 
  & \quad\quad =\big[ a\, (\varphi_{(1)}\re b \li
  \hat{S}^{-1}(\varphi_{(3)}))\,\Om^3  
       \big] \tp  \big[(\Om^2  \arr\varphi_{(2)}\arl
       \Om^4)(\Om^1\arr\psi\arl\Om^5)  \big].
\end{align*}
where the second line holds if $\dH$ has a coalgebra structure dual to the algebra $H$,
in particular if $H$ is finite dimensional.
\end{definition}
In \cite{HN1} it has been shown that $A \reli \dH$ is a unital algebra
containing $A \equiv A \tp \hat{\e} $ as a unital subalgebra.%
 \footnote{ The proof in \cite{HN1} straightforwardly generalizes to infinite dimensions, in
which case the generating matrix $\RR$ above has to be viewed as a map $\RR:\dH\to A\reli \dH$.}
Defining the
{\em generating matrix} $\RR \in H \tp (A\reli \dH)$ by $\RR : = \sum_i b_i
\tp (\e\tp b^i)$, where $\{b_i\}$ is a basis in $H$ with dual basis 
$\{b^i\}\subset \dH$, then $\RR$ obeys $(\ep\tp \id)(\RR) = \e_A\tp \hat{\e}$ and
\begin{align}
\label{4.7}
\RR\,[\e\tp  a] &= [a\1\tp a\0 ]\, \RR\, [\Si(a_{(-1)})\tp
\e_A], \\
\label{4.8}
 \RR^{13}\,\RR^{23} & = [\Om^4 \tp \Om^5\tp \Om^3 ]\, (\cop\tp\id)(\RR)
\, [\Om^2\tp\Om^1\tp \e_A].
\end{align}
Moreover, $A\reli \dH$ is the unique smallest unital algebra extension 
$B\supset A \equiv A \tp\hat{\e}$ such that there exists $\RR\in H\tp B$ 
satisfying $(\ep\tp\id)(\RR) = \e$ together with
\eqref{4.7} and \eqref{4.8}, see \cite{HN1}.

To prove the following Theorem we rely on the existence of a normalized cocentral biinvariant
form $\Sigma:H\tp H\to\kk$, which is why by Corollary~\ref{cor6.5} we require $H$ unimodular.
However, we don't know whether these conditions are also necessary.
\begin{theorem} 
  \label{thm4.8}
Let $H$ be a f.d. unimodular quasi-Hopf algebra and let   $(A,\delta,
\Psi)$ be a twosided $H$-comodule algebra. 
If  $A$ is semisimple Artinian and if $H$ admits a normalized left cointegral $\lambda_0\in \L$,
 then the
  diagonal crossed product $A \reli H$ is semisimple Artinian.
\end{theorem}
\begin{proof}
%Recall that a finite dimensional algebra is semisimple  if and only if all
%its finite dimensional left modules are reducible (Wedderburn's
%structure Theorem, see e.g. \cite{Abe}). Thus it suffices to 
We show that
for any two $(A\reli \dH)$-modules $V,W$, where $W \subset V$ is a
submodule, there exists a $(A\reli\dH)$-linear surjection $\bar{p}:V\to W$. 
Denoting the canonical embedding $i: W \hookrightarrow V$
and $\RR_W = (\id\tp\pi_W)(\RR) \in  H\tp\End (W)$ and $\RR_V =
(\id\tp\pi_V)(\RR)$
it therefore suffices to find a surjection $\bar{p}:V\to W$ satisfying
\begin{align}
  \label{4.9}
   \bar{p} \circ i &= \id_W, \\
  \label{4.10}
    \bar{p} \circ a &= a \circ \bar{p}, \quad \forall a \in A, \\
  \label{4.11}
    (\id\tp \bar{p})\circ \RR_V &= \RR_W \circ (\id\tp\bar{p}),
\end{align}
where at the l.h.s. and r.h.s. of \eqref{4.10} we have used the shortcut
notation $a \equiv \pi_V(a) \in \End (V)$ and $a \equiv \pi_W(a) \in \End
(W)$, respectively. This notation will also be used frequently below. We
proceed as follows.  Viewing $V$ and $W$ as $A$-modules, the semisimplicity
of $A$ implies the existence of an $A$-linear surjection $p: V \to W$,
satisfying \eqref{4.9} and \eqref{4.10}. Denoting $\Omb := \Om^{-1}$, we now
define the map $\bar{p}: V\to W$ in terms of $p$ by
\begin{equation*}
 \bar{p} : = (\Sigma_0\tp\id_W) \circ (\Omb^4 \tp \Omb^5 \tp \Omb^3)
    \circ \RR^{13}_W \circ (\id\tp\id\tp p)\circ \RR^{23}_V \circ
    (\Omb^2 \tp \Omb^1 \tp \id_V)\,,
\end{equation*}
where $\Sigma_0$ is the normalized biinvariant cocentral form associated with $\lambda_0$.
To show \eqref{4.9} note that $\RR_V \circ (\e\tp i) = (\e\tp
i)\circ \RR_W$, since by assumption the embedding $i$ is 
$A\reli\dH$-linear. Using $p\circ i = \id_W$ this implies
\begin{equation*}
  \RR^{13}_W \circ (\id\tp\id\tp p)\circ \RR^{23}_V \circ (\id\tp\id\tp
  i) = \RR^{13}_W \circ \RR^{23}_W
\end{equation*}
and therefore
\begin{align*}
  \bar{p} \circ i &= (\Sigma_0\tp \id_W) \circ (\Omb^4 \tp \Omb^5
    \tp \Omb^3) 
    \circ \RR^{13}_W \circ \RR^{23}_W \circ (\Omb^2 \tp \Omb^1 \tp
    \id_W) \\
    & = (\Sigma_0 \tp \id_W) \circ (\cop\tp\id)(\RR_W)\\
     & = (\ep\tp \id_W)(\RR_W)\\
     & = \id_W,
\end{align*}
where we have used \eqref{4.8}, then the normalization condition $\Sigma_0\circ\Delta=\ep$ 
and finally the identity $(\ep\tp\id)(\RR) = \e\tp \hat{\e}$.
Thus we have proven \eqref{4.9}.

$A$-linearity \eqref{4.10} follows from \eqref{4.7}
 and biinvariance of $\Sigma_0$, since one computes
 \begin{align*}
   \bar{p}\circ a & = (\Sigma_0 \tp \id_W) \circ\big(\cdots) \circ (\e\tp\e\tp
                  a) \\ 
                  & = \Big(\big(\cop(\Si(a_{(-1)}))\arr \Sigma_0 \arl
                  \cop(a\1)\big) \tp a\0 \Big) \circ (\cdots) \\
                 & =a\circ \bar{p}.
 \end{align*}
We are left to show \eqref{4.11}. For the l.h.s. we get 
\begin{align}
   \notag
  (\id\tp \bar{p})\circ \RR_V & =
   (\Sigma_0 \tp\id\tp\id_W) \circ (\Omb^4 \tp \Omb^5 \tp \id \tp
  \Omb^3) 
  \circ \RR^{14}_W \circ  \\
   &\quad \notag  (\id^3 \tp p) \circ
   \RR^{24}_V \circ \RR^{34}_V \circ
    (\Omb^2\tp\Omb^1 \tp \id\tp\id_V) \\
   &= \notag
    (\Sigma_0 \tp\id\tp\id_W) \circ 
   \big( \Omb^4 \Om^3\1 \tp \Omb^5\Om^4 \tp \Om^5 \tp
  \Omb^3\Om^3\0\big) \circ \RR^{14}_W \circ  \\
&\quad\notag   (\id^3 \tp p)\circ
  \big[(\cop\tp\id)(\RR_V)\big]^{234} \circ
    \big(\Si(\Om^3_{(-1)})\Omb^2 \tp \Om^2\Omb^1 \tp \Om^1 \tp
  \id_V\big) \\
  & = \notag
   (\Sigma_0 \tp\id\tp\id_W)\circ (\phi^{-1}\tp\id_W) \circ 
    \big( \Omb^4 \tp \cop(\Omb^5)\tp\Omb^3 \big) \circ    
    \RR^{14}_W \circ  \\
&\quad \label{4.12} 
  (\id^3 \tp p)\circ \big[(\cop\tp\id)(\RR_V)\big]^{234}
     \circ \big(\Omb^2 \tp \cop(\Omb^1)\tp\id_V\big) 
     \circ (\phi\tp\id_V). 
\end{align}
Here we have used \eqref{4.7},\eqref{4.8} and the identity 
\begin{align*}
  &\Om^1 \tp \Om^2 \Omb^1 \tp \Si(\Om^3_{(-1)}) \Omb^2 \tp \Omb^3
  \Om^3\0 \tp \Omb^4\Om^3\1 \tp \Omb^5\Om^4 \tp \Om^5 \\
 &\quad\quad =\Big[ \big( \cop^{op}(\Omb^1)\tp \Omb^2 \big)\phi^{321}
  \big(\Om^1\tp \cop^{op}(\Om^2)\big)\Big] \tp \Om^3\Omb^3 \\ 
  &\quad\quad\quad \tp
   \Big[\big(\cop(\Om^4) \tp \Om^5 \big) \phi^{-1} \big(\Omb^4 \tp
  \cop(\Omb^5)\big)\Big]
\end{align*}
- following from \eqref{4.6}, \eqref{4.5b} and \eqref{1.18} - and the
 biinvariance of $\Sigma_0$.
 A similar calculation yields for the r.h.s. of \eqref{4.11}
\begin{align}
  \label{4.13} & \RR_W \circ (\id\tp \bar{p}) =
      (\id\tp\Sigma_0 \tp \id_W) \circ (\phi\tp\id_W)
       \circ \big( \cop(\Omb^4) \tp \Omb^5 \tp \Omb^3\big) 
       \\
\notag
&\quad\quad \circ \big[ (\cop\tp\id)(\RR_W)\big]^{124}  \circ (\id^3 \tp
      \bar{p}) \circ \RR^{34}_V \circ
     \big( \cop(\Omb^2) \tp \Omb^1 \tp \id_V \big)\circ (\phi^{-1}\tp\id_V)
\end{align}
Comparing \eqref{4.12} and \eqref{4.13}, both expressions coincide
since $\Sigma_0$ is cocentral. This proves \eqref{4.11} and therefore
concludes the proof of Theorem~\ref{thm4.8}.
\end{proof}
Theorem~\ref{thm4.8} and Theorem~\ref{thm6.2} imply
\begin{corollary}
\label{cor4.9}
Let $H$ be a f.d. semisimple unimodular quasi-Hopf algebra, admitting
a normalized cointegral $\lambda_0$, i.e. $\lambda_0 (\beta S(\alpha)) = 1$. Then
the quantum double $D(H) = H \reli \hat{H}$ of $H$ is semisimple.
\end{corollary}

\noindent
{\bf Conjecture:}
More generally,
   for $\lambda\in\L$ and $r\in R$, we expect $(\beta\arr\lambda \reli r)\in
 D(H)$ to be a left integral, at least if $H$ is unimodular
\footnote{Here we have used the convention $D(H) = \dH \reli H$, see
\cite{HN1}}.  
Moreover, we expect $D(H)$ to be always unimodular as in the Hopf algebra case.
 Note that the counit on $D(H)$ is given by $\ep_D(\psi\reli a)=\psi(S^{-1}(\alpha))\ep(a)$.
Hence, under the conditions of Corollary~\ref{cor4.9} the  Haar integral in $D(H)$
should be given by $(\beta\arr\lambda_0 \reli e)$.

\bibliography{integral}
\bibliographystyle{amsalpha}

\end{document}